\documentclass[titlepage]{article}
\usepackage{amsmath}
\begin{document}
\title{A Surgery Theory for Manifolds of Bounded Geometry}
\author{Oliver Attie}
\maketitle
\section{Introduction}

In this paper we introduce a new application of geometric topology to the
study of Riemannian metrics. The purpose of this application is to classify
metrics of bounded geometry up to smooth quasi-isometry on an open manifold.

A manifold of bounded geometry is a non-compact manifold whose geometric
complexity is bounded. Such manifolds can be described metrically as having
sectional curvature bounded in absolute value and injectivity radius bounded
below. Cheeger's finiteness theorem is equivalent to saying that in the PL
sense such a manifold has a triangulation with a uniform bound on the number
of simplices in the link of each vertex. Universal coverings of compact
manifolds and leaves of foliations lie within this class of open manifolds.
In fact, Gromov has remarked that every manifold of bounded geometry is the
leaf of a lamination of the infinite dimensional compact space of Riemannian
metrics. It is still an open question, however, whether every manifold of
bounded geometry is the leaf of a foliation of a compact manifold. (If the
foliation is $C^1$ the answer is no, see \cite{attiehurder}).

Two non-compact manifolds are said to be smoothly quasi-isometric when there
exists a diffeomorphism $f$ between them so that their distance-metrics
satisfy
$$\frac{1}{c}d(x,y)\le d^\prime(f(x),f(y))\le cd(x,y).$$
The classification problem for open manifolds, analogous to the diffeomorphism
classification problem for compact manifolds is the following: Given $M$ and
$N$ and a suitable map $f:M \to N$, when is $f$ boundedly homotopic to a
smooth quasi-isometry?

The ``answer'' to the classification problem here is not an explicit one, but
rather one in the form of a surgery exact sequence, none of whose terms can
be computed in general, as is the case with compact manifolds. However, we
will be able to present the computation for a large class of examples. We
next describe the surgery theory (and correspondingly new notions of algebraic
topology) which will be developed here.

Surgery theory for high-dimensional manifold classification was introduced
by Milnor \cite{milnor1,milnor2,km} in the late 1950s and early 1960s.
The surgery exact sequence was first proven
by Kervaire and Milnor for homotopy spheres in \cite{km},
and for high-dimensional manifolds by Browder \cite{browder1,browder2},
Novikov \cite{novikov1,novikov2}, Sullivan \cite{ranicki2}, 
Casson \cite{ranicki2}
and Wall \cite{wall}. Freedman \cite{freedman1,freedman2,fq} extended surgery
classification to 4-manifolds with some restrictions on the fundamental
group. 
Connell and Hollingsworth
\cite{connellhollingsworth} introduced controlled algebraic topology
in the 1960s. Anderson and Hsiang \cite{ah}, Chapman \cite{chapman1,chapman2}, 
Ferry \cite{ferry1}, Quinn \cite{quinn1,quinn2,quinn3}, 
Pedersen and Weibel \cite{pw} developed $\epsilon$-controlled and boundedly  
controlled topology in the 1970s and
1980s. Ferry and Pedersen \cite{fp} introduced boundedly controlled
surgery in the late 1980s. Our result on the uniqueness for $\textbf{R}^n$
extends the result proven by Siebenmann in 1968 \cite{siebenmann1}.

The surgery theory developed in this paper is an L-theoretic analogue of the
index theorem of Roe \cite{roe1} in the sense that both Roe's index and 
the surgery obstruction lie in groups that, if a Baum-Connes \cite{baumconnes}
 or Borel type 
conjecture were true (Baum-Connes for Roe's coarse theory has been shown
to be false by Higson, Lafforgue and Skandalis \cite{hls}),
could be expressed as $L^\infty$ homology with 
coefficients in a spectrum. Weinberger observed in 1990 that boundedly 
controlled surgery should be analogous to Roe's coarse index theory  
\cite{roe2}.
Roe's coarse index theory \cite{roe2, higsonroe} 
has been shown to be related to boundedly controlled surgery for some
spaces in the sense that both the index and the surgery obstruction lie
in exotic homology with coefficients in a spectrum \cite{dfw}.
A character map (following Connes and Moscovici \cite{connesmoscovici})
 from the cyclic homology of Roe's uniformly smoothing algebra to 
uniformly finite homology has been constructed by Block and Weinberger
\cite{bw2}.Roe's coarse index theory has also been used to prove the 
homotopy invariance of rational Pontrjagin classes (originally due to 
Novikov \cite{novikov3}), see \cite{prw}.
 
The main results of this classification theory for $bg$ manifolds contrast
sharply with other known results on the topology of non-compact manifolds.
For example, Euclidean and hyperbolic space are homeomorphic to each other
in the boundedly controlled category. However, the quasi-isometric 
classification
of universal covers of manifolds with fundamental group a surface group
exhibits the following phenomenon: there exist manifolds $X$ and $Y$ and a
homotopy equivalence $f:X \to Y$ so that $f^*(p_i)-p_i \neq 0$, where $p_i$
denotes the $i$-th Pontrjagin class, and $f$ lifts to a map which is boundedly,
but non-equivariantly, homotopic to a quasi-isometry on the
universal covers. In distinction to this, every homotopy equivalence of 
manifolds
with free abelian fundamental group which lifts to a map which is 
boundedly homotopic to a  quasi-isometry must
preserve Pontrjagin classes; see \cite{abw}.

In \cite{attie} we introduced a new PL category, with objects simplicial
complexes with bounded combinatorial complexity, and maps with bounded combinatorial
complexity (see Definitions 2.1-2.14). In this category, hyperbolic $n$-space $\textbf{H}^n$ has a
different ``homotopy type" (see Definitions 2.15-2.17) than euclidean $n$-space $\textbf{R}^n$. There
is an invariant of this ``homotopy type", the uniformly finite homology,
denoted $H^{uff}_*(X;G)$ which is defined in Section 3, where $G$ is an abelian group equipped
with a norm, so that
$$H^{uff}_0(\textbf{H}^n;\textbf{Z})=0$$
and
$$H^{uff}_0(\textbf{R}^n;\textbf{Z}) \neq 0.$$
This homology theory can be thought of as $L^\infty$-homology with coefficients
in $G$. In addition to this notion of a ``homotopy type" we introduced the
notion of a ``simple homotopy type", discussed in section 4, with a smooth quasi-isometry being an
example of this new type of ``simple homotopy equivalence". We also introduce
the notion of a ``structure set" (see Definition 5.10) in this new category,
$\mathcal{S}^{bg,s}_{TOP}(X)$, which is the set of ``homotopy equivalences" of
manifolds of bounded geometry to $X$ in this category, modulo ``homeomorphisms"
(see Definition 2.11) in this category.

The main classification result proven here is:
\newtheorem{theorem}{Theorem}[section]
\begin{theorem}
 Let $M^k$ be a compact manifold. Then the bg simple structure set
of $M \times R^n$, $k+n \ge 5$ is
$$\mathcal{S}^{bg,s}_{TOP}(M \times \textbf{R}^n)=H^{uff}_0(\textbf{R}^n;\mathcal{S}_{TOP}(M \times D^n,\partial))
\oplus ... \oplus H^{uff}_n(\textbf{R}^n;\mathcal{S}_{TOP}^{2-n}(M))$$
where $\mathcal{S}_{TOP}^{-i}(M)$ denotes the fiber of the assembly map of Ranicki's
lower L-theory, and we take the convention $\mathcal{S}_{TOP}^1(N,\partial N)=
\mathcal{S}_{TOP}^h(N,\partial N)$, where $N$ is a compact PL manifold with
boundary $\partial N$.
\end{theorem}

From this we derive two results:

\begin{theorem}
Let $M^n$, $n \ge 5$ be a uniformly contractible smooth manifold
of bounded geometry. Suppose further that there is a surjective map
$$f:M \to \textbf{R}^n$$
which is EPL in the sense of \cite{bw} or a coarse map of bounded
geometry in the sense of \cite{attie}. Then $M$ is 
smoothly quasi-isometric to $\textbf{R}^n$.
\end{theorem}

We also have the following result:

\begin{theorem}
Let $M^k$ be a compact manifold and
$$f:M \times T^n \to N^{k+n}$$
a homotopy equivalence, $k+n \ge 5$. Then the free abelian cover of $f$
$$\tilde{f}:M \times \textbf{R}^n \to \tilde{N}$$
is $bg$ homotopic to a quasi-isometry if and only if the lift of $f$ to a
finite cover is homotopic to a diffeomorphism.
\end{theorem}

I would like to thank Shmuel Weinberger for a number of fundamental discussions
and suggestions, which were crucial to the development of this work. I would
also like to thank Sylvain Cappell for discussions and for suggesting the
filtration of the Whitehead group discussed below. Finally, I would like to
thank Andrew Ranicki for pointing out that the results in this paper can be
obtained through the application of algebraic surgery methods and without
using topological manifolds or surgery spectra. I would also like to thank
Hans Munkholm for assistance in clarifying some of the definitions.

\section{Preliminaries}

In this section we discuss the general results needed for the quasi-isometry
classification of manifolds of bounded geometry. Bounded geometry was
first studied by Cheeger and Gromov in \cite{cheegergromov}.
We recall the definitions of
simplicial complexes of bounded geometry, homotopy equivalences of bounded
geometry and notions of controlled topology which will be used in the paper.
These definitions are taken from \cite{attie}.
\newtheorem{definition}{Definition}[section]
\begin{definition}
A simplicial complex $X$ has bounded geometry if there is a
uniform bound on the number of simplices in the link of each vertex of $X$.
\end{definition}

\begin{definition}
A simplicial map $f:X \to Y$ of simplicial complexes of
bounded geometry is said to have bounded geometry if the inverse image of each
simplex $\Delta$ of $Y$ contains a uniformly bounded number of simplices
of $X$. The uniform bound is called the complexity of the map.
\end{definition}

For continuous maps, there is a notion of bounded geometry, which can be
found in \cite{bw} where it is called EPL:

\begin{definition}
Let $X$ and $Y$ be metric spaces. A coarse map of bounded
geometry is (not necessarily continuous) map $f:X \to Y$ satisfy the conditions:

i. A condition similar to a uniform Lipschitz condition. That is, given $r>0$,
there is a uniform $s>0$ depending only on $r$ so that $f(B(x,r)) \subset
B(f(x),s)$, where $B(x,r)$ denotes the metric ball of radius $r$ around $x$.

ii. It is effectively proper. That is, given $r>0$ there exists a uniform
$s>0$ depending only on $r$ so that $f^{-1}(B(f(y),r) \subset B(y,s).$
\end{definition}

We next recall the conditions for bounded geometry on the Riemannian metric
of a smooth manifold.

\begin{definition}
A complete Riemannian manifold $M$ is said to have bounded
geometry if its injectivity radius $inj_M > c > 0$ for some constant $c$
and its sectional curvature is bounded in absolute value. Recall that the
injectivity radius of a complete Riemannian manifold is the infimum of the
injectivity radii at each point of $M$. The injectivity radius at a point is
the maximum radius for which the exponential map is injective.
\end{definition}

\begin{definition}
A smooth map of bounded geometry is a smooth map which is
effectively proper so that the $C^2$ norm of $f$ is uniformly bounded.
\end{definition}

\begin{definition}
A subdivision of a simplicial complex of bounded geometry is
said to be uniform if

i. Each simplex is subdivided a uniformly bounded number of times on its
$n$-skeleton, where the $n$-skeleton is the union of $n$-dimensional
sub-simplices of the simplex.

ii. The distortion
$$sup(length(e), length(e)^{-1})$$
of each edge $e$ of the of the subdivided complex is uniformly bounded in the
metric given by the barycentric coordinates of the original complex.
\end{definition}

\begin{definition}
A metric space $P$ is a $bg$ polyhedron if:

i. It is topologically a subset $P \subset R^n$.

ii. Each point $ a \in P$ has a cone neighborhood $N=aL$ of $P$ in the given
Euclidean space, where $L$ is compact and there is a uniform upper bound
for all $a \in P$ for the number of simplices needed to triangulate $L$.
\end{definition}

\begin{definition}
A map $f:P \to Q$ between $bg$ polyhedra is $bg$ PL if it
is piecewise linear and has bounded distortion, i.e. the distortion of the image
of a simplex is uniformly bounded. This is equivalent to saying that the
graph of $f$ is a $bg$ polyhedron.
\end{definition}

\begin{definition}
A PL manifold of bounded geometry is a $bg$ polyhedron so
that each point $x \in M$ has a neighborhood in $M$ which is PL homeomorphic
to an open set of $R^n$, with a uniform bound on the distortion of the PL
homeomorphism over $M$.
\end{definition}
\newtheorem{remark}{Remark}[section]
\begin{remark} A PL map of bounded geometry is an equivalence class of
simplicial maps of bounded geometry under the equivalence of uniform subdivision.
This follows by writing $bg$ polyhedra as unions of simplices.
\end{remark}
\begin{definition}
Let $f:M \to N$ be a smooth map between Riemannian manifolds.
Then $f$ has bounded dilatation if there is a constant $C$ so that
$$\mid f_* v \mid \le C \mid v \mid$$
for all $v \in TM$. $f$ is called a smooth
quasi-isometry if it is a diffeomorphism and both $f$ and $f^{-1}$ have
bounded dilatation.
\end{definition}

We shall also need to use the following:
\begin{definition}
Let $f:M \to N$ be a continuous map between PL manifolds of bounded geometry,
then $f$ is a $bg$ homeomorphism (or equivalently a continuous quasi-isometry)
if $f$ has a continuous inverse $f^{-1}$
such that the distance metrics, $d$ of $M$ and $d^\prime$ of $N$ satisfy
$$\frac{1}{c}d(x,y) \le d^\prime(f(x),f(y)) \le cd(x,y)$$
\end{definition}

We recall the following Theorem due essentially to Cheeger, M\"uller and
Schrader \cite{cheeger, cms}, from \cite{attie}:

\begin{theorem}
Let $M$ be a smooth manifold with a Riemannian metric of \break bounded
geometry. Then $M$ admits a triangulation as a simplicial complex of bounded
geometry whose metric given by barycentric coordinates is quasi-isometric to
the metric on $M$ induced by the Riemannian structure. This triangulation is
unique up to uniform subdivision. Conversely, if $M$ is a simplicial complex
of bounded geometry which is a triangulation of a smooth manifold, then this
smooth manifold admits a metric of bounded geometry with respect to which it
is quasi-isometric to $M$.
\end{theorem}
\newtheorem{corollary}{Corollary}[section]
\begin{corollary}
A smooth map which can be simplicially approximated by a
simplicial map of bounded geometry for appropriate triangulation of the source
and target, can be approximated by a smooth map of bounded geometry. Conversely,
any smooth map of bounded geometry can be simplicially approximated by a PL
map of bounded geometry.
\end{corollary}
\begin{definition}
Let $M^n \subset N^{n+q}$, then $N$ is an abstract regular
neighborhood of bounded geometry if $N$ collapses via a $bg$ map to $M$.
\end{definition}
\begin{definition}
A bounded geometry $q$-block bundle $\xi^q$ consists of a
total space $E(\xi)$ and a $bg$ simplicial complex $K$ so that $\mid K \mid
\subset E(\xi)$ satisfying

i. For each $n$-cell $\sigma_i \in K$, there exists an $(n+q)$-ball $\beta_i
\subset E(\xi)$ so that $(\beta_i,\sigma_i) \simeq (I^{n+q},I^n)$

ii. $E(\xi)$ is the union of blocks $\beta_i$.

iii. The interiors of blocks are disjoint.

iv. They are compact polyhedra and fall into a finite number of types (as
simplicial complexes).

v. Let $L=\sigma_i \cap \sigma_j$, then $\beta_i \cap \beta_j$ is the bounded
union of blocks over cells of $L$.

$\xi^q,\eta^q/K$ are $bg$ isomorphic if there is a $bg$ homeomorphism $h:E(\xi)
\to E(\eta)$, $h \mid K=1$, $h(\beta_i(\xi))=\beta_i(\eta)$, $\sigma_i \in K$.
$\xi \sim \eta$ or $\xi$ equivalent to $\eta$, if there exist uniform subdivisions
$\xi^\prime, \eta^\prime$ so that $\xi^\prime \simeq \eta^\prime$.
Let $X=\mid K \mid$. Let
$I_q(K)$ denote the set of $bg$ isomorphism classes of $q$-block bundles over
$K$, $I_q(X)$ the set of $bg$ equivalence classes over $X$. Then amalgamation
gives a bijection between the two sets.
\end{definition}
The following Theorem is from \cite{attie}.
\begin{theorem}
Let $N^{n+q}$ be a $bg$ abstract regular neighborhood of $M^n$
and suppose $L \subset K$ so that $(M, \partial M)=(\mid K \mid, \mid L \mid)$.
Then there is $\xi^q/K$ with $E(\xi)=N$.
\end{theorem}
\begin{definition}
Let $M,N \subset Q$ be $bg$ submanifolds of the $bg$ manifold
$Q$, all in the $bg$ PL category. Let $\xi$ a normal $bg$ block bundle on $M$.
Then $N$ is $bg$ transverse to $M$ with respect $\xi$ if there is a uniform
subdivision $\xi^\prime$ of $\xi$ so that $N \cap E(\xi)=E(\xi^\prime \mid
N \cap M)$.
\end{definition}

\begin{theorem}[Transversality Theorem \cite{attie}] Let $M,N \subset Q$ be $bg$ submanifolds of the
$bg$ manifold $Q$. There is an ambient isotopy of bounded geometry of $Q$
carrying $N$ $bg$ transverse to $M$ with respect $\xi$.
\end{theorem}
\textit{Proof.} This is word for word as in \cite{rourkesanderson}, with
bounded geometry in front of every term except subdivision, which must have
the word uniform in front of it instead. We recall that the strategy of proof
in \cite{rourkesanderson} is to use the Zeeman unknotting theorem to induct
on the skeleta of the dual triangulation of $M$. The analogous unknotting
theorem for the bounded geometry case states that if one has an infinite
collection of $bg$ embedded spheres, each of which can be unknotted in the
PL category, then one can unknot them by a $bg$ PL isotopy. The rest also
goes through verbatim as stated above.

\bigskip\noindent
\textbf{Example:} Consider the submanifold $y=e^{-\mid x \mid}sin(x)$ of $\textbf{R}^2$, which
is a smooth submanifold of bounded geometry. This submanifold intersects the 
line $y=0$ transversely in an infinite number of points, but is not $bg$ 
transverse. The embedding of this submanifold into $\textbf{R}^2$ cannot be 
simplicially as a $bg$ PL embedding so that the intersection is still a 
countably infinite number of points, since this requires $\textbf{R}^2$ to be
triangulated in such a way that the volume of a simplex is not bounded below.
However this embedding can still be $C^0$-approximated by a $bg$ PL
map which is $y=0$ outside a compact set, which can be made arbitrarily large. 
This is clearly not $bg$ transverse since the dimension of the intersection
is positive outside the compact set. If we isotop this curve by adding a small
constant to $y$ in the equation above, it becomes $bg$ transverse.
\begin{definition} A homotopy of bounded geometry between two maps $f_0$ and 
$f_1$ of bounded geometry between simplicial complexes $X$ and $Y$ of bounded
geometry is a map of bounded geometry $F: X \times I \to Y$ so that
$F \mid X \times 0=f_0$ and $F \mid X \times 1 =f_1$. We write this
$f_0 \sim_{bg} f_1$. 
\end{definition}
\begin{definition} A homotopy equivalence of bounded geometry is a map $f$ of
bounded geometry so that there is a map $g$ of bounded geometry with $f \circ
g$ and $g \circ f$ $bg$ homotopic to the identity.
\end{definition}

\begin{definition} A CW-complex of bounded geometry is defined to be a CW-complex
with a uniformly bounded number of cells attached to each cell and a finite
number of homeomorphism types of attaching maps. A $bg$ $n$-cell is a discrete
collection of $n$-cells $\Sigma \times I^n$, equipped with an attaching map
$\psi:\Sigma \times I^n \to X$. Two attaching maps $\psi_1,\psi_2: \Sigma
\times I^n \to X$ are of the same homeomorphism type if there is a cellular
homeomorphism $h:X \to X$ so that $h \psi_1 h^{-1}=\psi_2$.
\end{definition}

\begin{definition} Let $X_1$ and $X_2$ be spaces equipped with continuous maps
$p_1,p_2$ to a metric space $Z$. Then a map $f:X_1 \to X_2$ is boundedly
controlled if there exists an integer $m \geq 0$ so that for all $z \in Z$,
$r \geq 0$, $p_1^{-1}(B_r(z)) \subseteq f(p_2^{-1}(B_{r+m}(z))$, where
$B_r(z)$ denotes the metric ball in $Z$ of radius $r$ about $z$. Or,
equivalently, there is a constant $m \geq 0$ so that
$$dist_Z(p_2 \circ f(x), p_1(x))< m$$
for all $x \in X$.
\end{definition}
The next proposition is from \cite{attie}.
\newtheorem{proposition}{Proposition}[section]
\begin{proposition}
 Let $X$ and $Y$ be a simplicial complexes of bounded geometry
equipped with a map of bounded geometry to a simplicial complex $Z$ of
bounded geometry. Then any map $f:X \to Y$ is boundedly controlled only if it
has bounded geometry.
\end{proposition}
\textit{Proof.} The property of having bounded geometry for a map is similar to
injectivity, except that the inverse image is uniformly bounded, rather
than being a point. So if $q \circ f$ is of bounded geometry, $f$ must be.
Hence if $f$ is boundedly controlled, then $\mid q \circ f(x) -p(x) \mid <M$
for all $x$, so that by the bounded geometry of the triangulation, $q \circ f$
has bounded geometry, so that $f$ does also.

\bigskip\noindent
We shall use this result also for smooth maps, but in this case one must take
care, as the following example shows. Consider the map $f(x)=x^2$ and the
map $g(x)=\sqrt{x}$ on $\textbf{R}$. Neither of these maps is $bg$. However,
$f \circ g$ is $bg$. This seems to contradict both the observation above
and the correspondance between smooth and PL maps. However, this example fails
because neither $f$ nor $g$ are simplicial with respect to a $bg$ triangulation
of $\textbf{R}$. Since smooth control maps will be used in the sequel, this
phenomenon must be taken into account and smooth maps will be required to be
simplicial with respect to a $bg$ triangulation.

We note also that the converse of this proposition is false if taken in its
most literal sense: bounded geometry does not imply bounded control. For
example, consider multiplication by 2 on $\textbf{R}$, controlled over itself
by the identity. However, every map of bounded geometry $f:X \to Y$ can be
considered to be boundedly controlled over $Y$, taking the control maps to be
$f$ and the identity, respectively. We will thus consider the two notions as
equivalent from now on. 

The following definitions are due to Anderson and Munkholm \cite{am}.

\begin{definition}
Let $X$ be a space controlled over a metric space $Z$ by
a control map $p$. Denote by $\mathcal{P}$ the category of metric balls in $Z$ with
morphisms given by the inclusions. Define $\mathcal{P} G_1(X)$ to be the
category whose objects are pairs $(x,K)$ where $K \in \mid \mathcal{P} \mid$ is
an object of $\mathcal{P}$ and a morphism $(x,K) \to (y,L)$ is a pair $(\omega,i)$
where $i \in \mathcal{P}(K,L)$ is a morphism in $\mathcal{P}$ from $K$ to $L$ and
$\omega$ is a homotopy class of paths in $p^{-1}(L)$ from $y$ to $p^{-1}(i(x))$.
\end{definition}

\begin{definition}
The controlled homotopy groups $\pi^c_n(X,p)$ are defined
to be the functor
$$\pi_n^c(X):\mathcal{P} G_1(X) \to \mathcal{C}$$
where $\mathcal{C}$ is the category of pointed sets, groups or abelian groups
defined by setting
$$\pi_n^c(X,p)(x,K)=\pi_n(p^{-1}(K),x)$$
and $\pi_n^c(\omega,i)$ is the composite of the change of basepoint isomorphism
$\omega_*$ induced from $\omega$ and the homomorphism induced from the
inclusion $i$.
\end{definition}

\begin{definition}
The controlled homology $H^c_n(X,p)$ (with integer coefficients)
of a space controlled via $p:X \to Z$ is defined to be the pro-system \break
$H_n(p^{-1}(B(r,z))$ via the maps $B(r,z) \to B(r+1,z)$.
\end{definition}
\begin{definition} If $(X,p)$ and $(Y,q)$ are spaces controlled over $Z$ by
control maps $p$ and $q$ respectively, then $X$ is coextensive with $Y$ if
there exists an integer $m \geq 0$ so that if $p^{-1}(B_r(z)) \neq 0$ then
$q^{-1}(B_{r+m}(z)) \neq 0$ and the same with the roles of $p$ and $q$ reversed.
\end{definition}

We shall also use the category $\mathcal{C}_M(R)$ for a metric space $M$ 
introduced by Pedersen and Weibel \cite{pedersen1,pw}. 
It is shown by Anderson and 
Munkholm \cite{am},pp.263-4, that in the case where the metric space $M$
is path connected, the metric is proper, and satisfies the 
condition that if $B(r,z) \subset B(s,z)$ then $B(r+1,z) \subset B(s+1,z)$,
which is clearly satisfied for the case where $M$ is a manifold of bounded 
geometry, then $\mathcal{C}_M(\textbf{Z}\pi_1(X))$ is equivalent to finitely
generated free modules over $\textbf{Z}\mathcal{P}G_1(X)$, in the sense that
there is an additive functor between them which is an equivalence of categories,
which is natural with respect to boundedly controlled maps. This additive
functor induces an isomorphism on algebraic K-theory, which is natural with 
respect to boundedly controlled maps.

The next six definitions are from \cite{fp}.
We define controlled chains and cochains as in \cite{fp} as follows:
\begin{definition} An object $A$ in $\mathcal{C}_M(R)$ is a collection of
finitely generated free right $R$-modules $A_x$, one for each $x \in M$,
such that for each ball $C \subset M$ of finite radius, only finitely many
$A_x$, $x \in C$ are nonzero. A morphism $\phi: A \to B$ is a collection of
automorphisms $\phi^x_y: A_x \to B_y$ such that there exists $k=k(\phi)$
such that $\phi^x_y=0$ for $d(x,y) > k$.

The composition of $\phi:A \to B$ and $\psi: B \to C$ is given by
$(\psi \circ \phi)^x_y=\sum_{z \in M} \psi_y^x \phi_z^x$. The composition
$(\psi \circ \phi)$ satisfies the local finiteness and boundedness conditions
whenever $\psi$ and $\phi$ do.
\end{definition}
\begin{definition} The dual of an object $A$ in $\mathcal{C}_M(R)$ is the object
$A^*$ with $(A^*)_x=A^*_x=Hom_R(A_x,R)$ for each $x \in M$. $A^*_x$ is
naturally a left $R$-module, which we convert to a right $R$-module by
means of the anti-involution. If $\phi: A \to B$ is a morphism, then
$\phi^*: B^* \to A^*$ and $(\phi^*)^x_y(h)=h \circ \phi^x_y$, where
$h:B_x \to R$ and $\phi^y_x: A_y \to B_x$, $\phi^*$ is bounded whenever
$\phi$ is. Again, $\phi^*$ is naturally a left module homomorphism which
induces a homomorphism of right modules $B^* \to A^*$ via the anti-involution.
\end{definition}
\begin{definition} Consider a map $X \to M$

i. The map $p:X \to M$ is eventually continuous if there exist $k$ and a
covering $\{U_\alpha\}$ of $X$, such that the diameter of $p(U_\alpha)$
is less than $k$.

ii. A bounded CW complex over $M$ is a pair $(X,p)$ consisting of a CW
complex $X$ and an eventually continuous map $p: X \to M$ such that there
exists $k$ such that $diam(p(C)) < k$ for each cell $C$ of $X$. $(X,p)$ is
called proper if the closure of $p^{-1}(D)$ is compact for each compact
$D \subset M$. We consider $(X,p_1)$ and $(X,p_2)$ to be the same, if there
exists $k$ so that $d(p_1(x),p_2(x))<k$ for all $x$.
\end{definition}

\begin{definition} Consider a bounded CW complex $(X,p)$

i. The bounded CW complex $(X,p)$ is (-1)-connected if there is a $k \in
\textbf{R}_+$ so that for each point $m \in M$, there is a point $x \in X$
such that $d(p(x),m) < k$.

ii. $(X,p)$ is 0-connected if for every $d>0$ there exist $k=k(d)$ so that if
$x,y \in X$ and $d(p(x),p(y)) \le d$, then $x$ and $y$ may be joined by a path
in $X$ whose image in $M$ has diameter $< k(d)$. Note that 0-connected
does not imply -1-connected.
\end{definition}

\begin{definition} Let $p:X \to M$ be 0-connected, but not necessarily
(-1)-connected.

i. $(X,p)$ has trivial bounded fundamental group if for each $d > 0$ there
exist $k=k(d)$ so that for every loop $\alpha: S^1 \to X$ with
$diam(p \circ \alpha(S^1)) < d$, there is a map $\overline{\alpha}: D^2 \to X$
so that the diameter of $p \circ \overline{\alpha}: D^2 \to X$ so that the
diameter of $p \circ \overline{\alpha}(D^2)$ is smaller than $k$.

ii. $(X,p)$ has bounded fundamental group $\pi$ if there is a $\pi$ cover
$\tilde{X}$ so that $\tilde{X} \to M$ has trivial bounded fundamental group.
\end{definition}

\begin{definition} If $X$ is a CW complex, we will denote the cellular chains
of $\tilde{X}$ by $C_\#(X)$ considered as a chain complex of free right
$\textbf{Z}\pi_1(X)$-modules. When $p:X \to M$ is a proper bounded CW complex
with bounded fundamental group, we can consider $C_\#(X)$ to be a chain
complex in $\mathcal{C}_M(\textbf{Z}\pi_1(X))$ as follows: For each cell
$C \in X$, choose a point $c \in C$ and let $D_\#(X)_y$ be the free submodule
of $C_\#(X)$ generated by cells for which $p(c)=y$. The boundary map is
bounded, since cells have a fixed maximal size. We will denote the cellular
chains of $\tilde{X}$ by $D_\#(X)$ when we consider them as a chain complex
in $\mathcal{C}_M(\textbf{Z}\pi_1(X))$ and by $C_\#(X)$ when we consider
them as an ordinary chain complex of $\textbf{Z}\pi_1(X)$ modules. We will
denote $D_\#(X)^*$ by $D^\#(X)$. If $(X, \partial X)$ is a bounded CW pair,
$D_\#(X,\partial X)$ denotes the relative cellular chain complex regarded
as a chain complex in $\mathcal{C}_M(\textbf{Z}\pi_1(X))$.
\end{definition}
The following two theorems are modifications of the corresponding ones in
\cite{am} and can be found in \cite{attie}.

\begin{theorem}[Whitehead Theorem]
If $f:(X,p) \to (Y,q)$ is a map of bounded geometry of
CW complexes of bounded geometry, controlled by maps to $Z$ of bounded
geometry, then $f$ is a $bg$ homotopy equivalence if $(Y,q)$ is coextensive
with $(X,p)$ and for all $n \geq 0$, $f_*:\pi^c_n(X,p) \to f^!\pi^c_n(Y,q)$
is an isomorphism.
\end{theorem}

\begin{theorem}[Hurewicz Theorem]
 If $X$ is a simplicial complex of bounded geometry, then

i. $\pi^c_1(X)^{ab}=H_1^c(X)$

ii. If $\pi^c_i(X)=0$ for $i \leq n-1$ and $n \geq 2$ then, $H_i^c(X)=0$ for
$i \leq n-1$, and $H^c_n(X)=\pi^c_n(X)$.
\end{theorem}

\section{Uniformly Finite Homology}

In this section we review $L^\infty$ cohomology and uniformly finite homology as
defined by Gromov, Block-Weinberger \cite{bw}, Roe \cite{roe1} and Gersten \cite{gersten1,gersten2}.
Whyte \cite{whyte} has shown that the Poincar\'e dual of the $\hat{A}$-class in
0-dimensional uniformly finite homology is an obstruction to the existence
of a metric of positive scalar curvature on an open manifold.
 
\begin{definition}
 A normed abelian group is an abelian group $G$ equipped with
a norm function:
$$\mid \cdot \mid: G \to \textbf{R}_+$$
which is not necessarily continuous, but so that the induced function on the
Cayley graph of $G$ is non-decreasing as one moves away from the identity
element in $G$.
\end{definition}

\begin{definition}
Let $X$ be a $bg$ simplicial complex. The $i$-dimensional fine uniformly
finite homology groups of $X$ with coefficients in the normed group
\break $(G, \mid \cdot \mid)$, denoted $H_i^{uff}(X;G,\mid \cdot \mid)$  
are defined to be
the homology groups of the complex of infinite simplicial chains whose
coefficients are in $l^\infty$ with respect to the norm $\mid \cdot \mid$
on $G$. Define the group of $q$-chains $C_q^{uff}(X;G,\mid \cdot \mid)$ to be
the group of formal sums of $q$-simplices in $X$, $c=\sum a_\sigma \sigma$ so
that there exists $K > 0$ depending on $c$ so that $\mid a_\sigma \mid
\leq K$ and the number of simplices $\sigma$ lying in a ball of given size
is uniformly bounded. The boundary is defined to be the linear extension of
the simplicial boundary.
\end{definition}
We shall also need uniformly finite cohomology with coefficients in the 
normed group $(G,\mid \cdot \mid)$.
\begin{definition} Let $X$ be a $bg$ simplicial complex. The $i$-dimensional
fine uniformly finite cohomology groups of $X$ with coefficients in the normed
group $(G, \mid \cdot \mid)$, denoted $H^i_{uff}(X;G, \mid \cdot \mid)$ are
defined to be the cohomology groups of the complex of infinite simplicial 
cochains $c \in Hom(C_q(X;G),G)$ which satisfy $\mid c(\sigma) \mid \le K$
for all simplices $\sigma \in X$ and fixed $K > 0$ depending on $c$. The
coboundary is defined to be the simplicial coboundary.
\end{definition}
\begin{definition}
Let $f:X \to Y$ be a simplicial map of bounded geometry, and let
$(G, \mid \cdot \mid)$ be a normed group. Then the induced map
$$f_*:H_i^{uff}(X;G,\mid \cdot \mid) \to H_i^{uff}(Y;G,\mid \cdot \mid)$$
is defined by $f_*([c])=[f \circ c]$, where $c:C \to X$, $C$ a simplicial
complex and $c$ a $bg$ simplicial map, represents a class in $H_i^{uff}(X;
G,\mid \cdot \mid)$. This is well-defined because $f$ commutes with the
boundary.
\end{definition}

We now set a convention for the coefficients which will be used throughout
the rest of the paper: the group $\textbf{R}$ will always be normed with the
absolute value, as will $\textbf{Z}$. The groups $Wh(\pi)$, 
$K_i(\textbf{Z}\pi)$, $\mathcal{S}^{PL}(M \times D^n, \partial)$,
$\mathcal{S}^{TOP}(M \times D^n, \partial)$ can be 
written for the groups $\pi$ considered in this paper, as the countable 
direct sum of finite abelian groups with a free abelian group. The terms in
the direct sum will be given the absolute value norm as subsets of 
$\textbf{C}$. Assuming this convention, we will drop the symbol 
$\mid \cdot \mid$ in uniformly finite homology from now on.

\begin{definition} We recall for completeness the definition of locally
finite homology. Let X be a simplicial complex which is also a metric 
space. Define $C^{lf}(X;G)$, $G$ an abelian group, to be the group of 
formal sums of simplices in $X$, $c=\sum a_\sigma \sigma$, $a_\sigma \in G$
so that the number of simplices $\sigma$ lying in a particular ball in $X$
is bounded. Not that $H^{uff}_*(X;G)=H^{lf}_*(X;G)$ for any finite abelian
group $G$. 
\end{definition}

\bigskip\noindent
We now give some calculations of $H_*^{uff}(X;\textbf{Z})$ which show:

i. It depends on the metric structure of $X$ and not just on its topology:
e.g. $\textbf{H}^n$ and $\textbf{R}^n$ have different uf homology, even 
though they are diffeomorphic. 

ii. It can be very large: for $X=\textbf{R}$, $H_0^{uff}(X;\textbf{Z})$ is 
an uncountably generated $\textbf{R}$-module.

iii. There are deep connections between properties of $H^{uff}_*(X;\textbf{Z})$ and infinite
group theory. For example amenability, Gromov hyperbolicity and weak forms of 
rigidity are determined by $H_*^{uff}(X;\textbf{Z})$.

We will recall only the statements of the results in \cite{attie}. For
sketches of proofs see \cite{attie}.

\begin{proposition} The 0-dimensional uniformly finite homology group of 
$\textbf{R}$ is
$$H_0^{uff}(\textbf{R};\textbf{Z})=H_0^{uff}(\textbf{R};\textbf{R})=
\frac{\{\phi:\textbf{Z} \to \textbf{Z} \mid \mbox{ }\parallel \delta \phi 
\parallel_\infty < \infty\}}{\{\phi:\textbf{Z} \to \textbf{Z} \mid \mbox{ }
\parallel \phi \parallel_\infty < \infty\}}$$
where $\delta \phi(n)=\phi(n)-\phi(n-1)$ and $\parallel \cdot \parallel_\infty$
is the $L^\infty$ norm.
\end{proposition}
\begin{theorem} Let $X$ be a simply connected symmetric space of non-positive
curvature of rank $r$. Then
$$H_i^{uff}(X;\textbf{Z})=0$$
for $i \le n-r-1$
$$H_i^{uff}(X;\textbf{Z}) \neq 0$$
for $i \ge n-r$. 
\end{theorem}
\begin{proposition} Let $X$ be a manifold with curvature pinched between
two negative constants. Then
$$H_i^{uff}(X; \textbf{Z})=0$$
for $i<n-1$. 
\end{proposition}
Gersten \cite{gersten2} has proven theorems relating Gromov hyperbolicity to
the vanishing of $H_{n-2}^{uff}(X;\textbf{Z})$. 
See \cite{gersten2} for details.

\bigskip \noindent
Whether or not an infinite group is amenable is determined by its uniformly
finite homology. 

In defining amenability, we will ultimately quote \cite{bw} where essentially
the following theorem is proven: a $bg$ simplicial complex $X$ is non-amenable
if and only if $H_0^{uff}(X; \textbf{Z})=0$. 

However, before we present this definition, we will review the classical 
definitions of amenability.

\begin{definition} An infinite discrete group $\Gamma$ is said to be amenable
if and only if it satisfies the two equivalent conditions:

i. There is a bounded linear functional $\mu: l^\infty(\Gamma) \to \textbf{R}$
with $inf_{g \in \Gamma}(f(g)) \le \mu(f) \le sup_{g \in \Gamma}(f(g))$ and 
for all $g \in \Gamma$, $\mu(g \cdot f)=\mu(f)$, where $g \cdot f(x)
=f(g^{-1}x)$.

ii. For every $k$ in the interval $(0,1)$ and arbitrary finite set of elements
$a_1,...,a_n$ in $\Gamma$ there is a finite subset $E$ of $\Gamma$, $E \neq 0$,
so that
$$\#(E \cap a_i \cdot E) \ge k \#(E)$$
for $i=1,...,n$.
\end{definition}
Condition (i) is due to Von Neumann, and was the first criterion for 
amenability to be introduced. Condition (ii) is due to F\"olner and is known
as the F\"olner condition. 

The following criterion for amenability applies to finite dimensional 
simplicial complexes, and is due to Brooks and Gromov.
\begin{definition} An n-dimensional simplicial complex of bounded geometry
is said to be amenable if there is a sequence of n-dimensional compact 
subcomplexes $X_i \subset X$ so that the $X_i$ exhaust $X$ and 
$Vol \partial X_i/Vol X_i \to 0$, where $Vol \partial X_i$ means the number
of $n-1$-simplices on the boundary of $X_i$. 
\end{definition}

This criterion is related to amenability of groups in the following way:

\begin{proposition} An infinite discrete group $\Gamma$ is amenable if and 
only if the universal cover of any compact manifold with fundamental 
group $\Gamma$ is amenable. 
\end{proposition}

We now come to the theorem of Gromov and Block and Weinberger \cite{bw}.
Gromov found a criterion for the amenability of an infinite covering using
differential forms. Block and Weinberger generalized this to arbitrary 
simplicial complexes of bounded geometry. Recall that $H_\beta^*(X;\textbf{R})$
is the bounded de Rham cohomology of $X$ with real coefficients, where $X$
is a smooth manifold of bounded geometry.

\begin{definition} Denote by $\Omega^p_\beta(M)$ the Banach space of p-forms
on a complete, oriented Riemannian manifold $M$ which are bounded in the norm
$$\parallel \alpha \parallel=sup\{\mid \alpha(x) \mid + \mid d\alpha(x) \mid:
x \in M \}$$
This gives rise to a complex $d_i:\Omega^i_\beta(M) \to \Omega^{i+1}_\beta(M)$.
The bounded de Rham groups are defined by
$$H_\beta^p(M)=[\mbox{Ker }d_p]/[\mbox{Im }d_{p-1}]$$
Note that we are not taking the closure of $\mbox{Im }d$ in this definition.
\end{definition}
\begin{theorem} A non-compact manifold of smooth bounded geometry $X$ is \break 
amenable if and only if
$$H_0^{uff}(X; \textbf{Z}) \neq 0$$
or equivalently,
$$H_\beta^n(X; \textbf{R}) \neq 0$$
\end{theorem}

We now adopt for the purpose of the next proposition, the following criterion
for amenability, due to Block-Weinberger.
\begin{proposition} Let $X=B\Gamma$ be the classifying space of $\Gamma$ 
considered as a simplicial complex, $\tilde X$ an amenable covering of $X$. 
Then the map $H_*(X;\textbf{R}) \to H_*^{uff}(\tilde{X};\textbf{R})$ given
by taking each cycle to its lift to $\tilde{X}$ is injective.
\end{proposition}
\textit{Proof.} The following proof was suggested to the author by J.Block
and \break S.Weinberger. We use the invariant mean to construct a left inverse.
The main step of the proof is to show that the complexes $C_*^{uff}(\tilde{X};
\textbf{Z})$ and $C_*(\Gamma;l^\infty(\Gamma))$ are isomorphic. Here $\Gamma$
acts on $l^\infty(\Gamma)$ via $\gamma \cdot f(x)=f(\gamma \cdot x)$, for
any $f:\Gamma \to \textbf{R}$ in $l^\infty(\Gamma)$ and 
$C_*(\Gamma;l^\infty(\Gamma))$ is the standard bar resolution with respect to
this action. The main step is just a matter of unraveling the various 
definitions. Once we have the desired isomorphism, one constructs the map
$C_*(\Gamma; l^\infty(\Gamma)) \to C_*(X;\textbf{R})$ induced by applying
the invariant mean to the coefficients. This gives the desired left inverse
and the proposition follows.

\bigskip\noindent
We have produced an injection with real coefficients. Use of the de Rham
theorem proven below will result in an injection with rational coefficients.

We will prove next the de Rham and Poincar\'e duality theorems due to the 
author, J. Block and S. Weinberger, which also is proven in \cite{ab}. 
We first need a refinement of the notion of uniform subdivision. 
\begin{definition} A regular uniform subdivision is defined in the 
following manner. Let $\sigma=[p_0,...,p_m]$ be a simplex in $\textbf{R}^k$,
$k \ge m$. The vertices of the stanard subdivision $S\sigma$ of $\sigma$ are the 
points
$$p_{ij}=\frac{1}{2}(p_i+p_j), i \ge j.$$
Define a partial ordering of the vertices of $S\sigma$ by setting
$$p_{ij} \le p_{kl},$$
if $i \ge k$,$j \ge l$. The simplices of $S\sigma$ are the increasing 
sequences of vertices with respect to the above ordering.

We define a regular uniform subdivision to be any uniform subdivision which is 
a sequence of standard subdivisions.
\end{definition}
Regular uniform subdvision of a $bg$ simplicial complex $K$ produces a 
$bg$ simplicial complex $K^\prime$ and induces a map $s:C^*_{uff}(K;G) \to
C^*_{uff}(K^\prime;G)$ which has bounded norm.
\begin{theorem}[de Rham theorem]Let $M$ be a smooth $n$-dimensional manifold
of bounded geometry. Then there is an isomorphism
$$H_{uff}^i(M;\textbf{R}) \simeq H^i_\beta(M)$$
\end{theorem}
\textit{Proof.} The idea is to use Whitney and de Rham maps, with the
trinagulation constructed by Theorem 2.1 serving to show that these maps
are well-defined. The proof here imitates the ones found in \cite{dodziuk}, \cite{whitney}.
Define the de Rham map
$$\int:\Omega^i_\beta(M) \to C^i_{uff}(M;\textbf{R})$$
by
$$\int (\omega) \cdot \sigma=\int_\sigma \omega.$$
The de Rham map commutes with regular uniform subdivision. 
To define the Whitney map, which is the chain homotopy inverse of the de
Rham map, let $c_\sigma$ be the cochain which assigns the value 1 to 
$\sigma$ and the value 0 to every other simplex.

For each point $q$ in $M$ write $q$ in barycentric coordinates as 
$q=\sum \nu_\alpha(q)q_\alpha$, where $q_\alpha$ runs over the vertices of 
the triangulation, and $\alpha$ is the corresponding index. For each 
$\alpha$, let $Q_\alpha$ let $Q^\prime_\alpha$ be subsets of $M$ so that
$Q_\alpha$ is the set of all $p$ with $\nu_\alpha(p) \ge \frac{1}{n+1}$;
$Q^\prime_\alpha$ is the set of all $q$ with $\nu_\alpha(q) \le \frac{1}{n+2}$.
Then $Q_\alpha \subset star(q_\alpha)$, $Int(Q^\prime_\alpha \supset 
M-star(q_\alpha)$. Let $\phi_\alpha^\prime(p)$ be a smooth non-negative
real function in $M$ which has all of its derivatives uniformly bounded over
$\alpha$ and so that each is positive in $Q_\alpha$ and zero in 
$Q_\alpha^\prime$. Construct the partition of unity
$$\phi_\alpha(p)=\frac{\phi_\alpha^\prime(p)}{\sum_\beta \phi_\beta^\prime(p)}$$
We will choose normalizations of the $\phi^\prime_\alpha$ below, so we will
assume $\phi^\prime_\alpha$ normalized suitably for the moment.
Take any $p \in M$; since $p$ has at most $n+1$ non-zero barycentric 
coordinates at least one of these, say $\nu_\beta(p)$ is $\ge \frac{1}{n+1}$.
Hence $p \in Q_\beta$, $\phi^\prime_\beta(p) > 0$ and $\phi_\alpha(p)$ is
defined for all $\alpha$. Define the Whitney map $W$ on each $c_\sigma$ by
$$W(c_\sigma)=r!\sum_{i=0}^r(-1)^i\phi_{\alpha_i}d\phi_{\alpha_0}\wedge ...
\wedge d\hat{\phi}_{\alpha_i}\wedge ... \wedge d\phi_{\alpha_r},$$
where $\sigma=q_{\alpha_0}...q_{\alpha_r}$, and $\phi_\alpha$ is defined 
above. We can then extend by linearity. Note that
$$supp\mbox{ }W(c_\sigma) \subseteq \sigma.$$
Furthermore, because the triangulation is uniform, the resulting form is \break 
bounded. In fact, the map $W$ defines a bounded map on chains in the norm
defined by taking the supremum of the coefficients of the chain.

We now choose normalizations so that the de Rham theorem will be true. Choose
normalizations of the $\phi^\prime_\alpha$ so that
$$\int W(c_{\sigma_\alpha})\sigma_\beta=\int_{\sigma_\beta} W(c_{\sigma_\alpha})
=\delta^\beta_\alpha,$$
since $supp\mbox{ }W(c_\alpha) \subset \sigma_\alpha$. Using this 
normalization, we have inductively for $\sigma^\prime$ a face of $\sigma$,
$$\int_\sigma W(c_\sigma)=\int_\sigma W(c_{d\sigma})=\int_{\partial \sigma}
W(c_{\sigma^\prime})=\int_{\sigma^\prime}W(c_{\sigma^\prime})=1.$$
The normalize all of the lower dimensional cochains. Define $\int^*$ and 
$W^*$ to be the maps induced on cohomology by $\int$ and $W$ respectively.

We claim that $\int^*$ and $W^*$ are inverses of each other. In fact, by
Stokes' theorem we can observe that $\int$ and $W$ are chain maps. That
$\int W =Id$ is proved by observing that the normalization conditions prove the
result for each $c(\sigma)$ and then extending by linearity. That $W^*\int^*=
Id$ follows from a series of easy estimates. Let $\omega$ be a closed 
differential form on $\Omega_\beta^*(M)$. Suppose the cohomology class
$[\int \omega]=0$ in $H_\beta^{uff}(M;\textbf{R})$. Since integration commutes
with regular uniform subdivision, we have $[\int^\prime \omega]=0$, where 
$\int^\prime$ denotes the integration map over a uniformly subdivided
triangulation. Fix $\epsilon_1 > 0.$ We claim that there exists a regular
uniform subdivision so that
$$\parallel \omega-W \cdot \int \omega \parallel < \epsilon_1$$
Where $\parallel \cdot \parallel$ is the norm on $\Omega_\beta^*(M)$.
Note that we have the estimate
$$\mid \omega(x)-W \cdot \int \omega \mid \le C \cdot diam(\tau) \cdot
sup_{x \in \tau} \mid \frac{\partial \omega}{\partial x}\mid$$
over each simplex $\tau$ of the triangulation, where $C$ can be chosen 
independently of $\omega, \tau$. One the obtains an estimate in terms of 
the mesh of the triangulation $h$
$$\parallel \omega - W \cdot \int \omega \parallel \le 
\sum_{\tau \in K} \int_\tau \mid \omega(x)-W \cdot \int \omega(x) \mid dV$$
$$\le 4C \cdot h^N m \parallel \omega \parallel$$
where $m$ denotes the multiplicity of the intersections of the $bg$ coordinate
charts of $M$. Since $\int \omega$ is a boundary, we can find a cochain $f$
so that
$$\parallel \int \omega - \delta f \parallel \le \epsilon_2$$
Then
$$\parallel \omega-dWf\parallel \le \parallel \omega-W \cdot \int \omega 
\parallel + \parallel W \parallel \cdot \parallel \int \omega - \delta f 
\parallel$$
$$\le \epsilon_1+\parallel W \parallel \cdot \epsilon_2$$
Since $\epsilon_2$ is chosen independently of $\epsilon_1$ and both can be made
arbitrarily small, we obtain that the cohomology class $[\omega]=0$. This
proves the de Rham theorem.

We shall also prove Poincar\'e duality, following the proof in \cite{gh}.
\begin{theorem}[Poincar\'e duality]
 If $M$ is a manifold of bounded geometry, there
is an isomorphism
$$H^{n-i}_{uff}(M;\textbf{Z})=H_i^{uff}(M;\textbf{Z})$$ 
\end{theorem}
To prove the theorem we need the following definitions from \cite{ranicki1}:
\begin{definition} The $bg$ simplicial complex $K$ is ordered, so that
for each simplex $\sigma \in K$ the set
$$K^*(\sigma)=\{\tau \in K \mid \tau > \sigma, \mid \tau \mid=\mid \sigma \mid
+1\}$$
$$K_*(\sigma)=\{\tau \in K \mid \tau < \sigma, \mid \tau \mid=\mid \sigma \mid
-1 \}$$
\end{definition}
\begin{definition} The star and link of a simplex $\sigma \in K$ in a 
simplicial complex $K$ are the subcomplexes defined by
$$star_K(\sigma)=\{\tau \in K \mid \sigma\tau \in K\}$$
$$link_K(\sigma)=\{\tau \in K \mid \sigma\tau \in K, \sigma \cap \tau=\emptyset
\}$$
where $\sigma\tau$ is the simplex spanned by $\sigma$ and $\tau$.
The dual cell of $\sigma$ is the contractible subcomplex of the barycentric
subdivision $K^\prime$ defined by
$$D(\sigma,K)=\{\hat{\sigma}_0\hat{\sigma}_1...\hat{\sigma}_p \in K^\prime
\mid \sigma \le \sigma_0 < \sigma_1 < ... < \sigma_p\}$$
where $\hat{\sigma}$ is the barycenter of $\sigma$. 
The barycentric subdivision of the link of $\sigma \in K$ is isomorphic to 
the boundary of the dual cell $D(\sigma,K)$
$$(link_K(\sigma))^\prime \simeq \partial D(\sigma,K)$$
The star and link in $K^\prime$ of the barycentre $\hat{\sigma} \in K^\prime$
of $\sigma \in K$ are given by the joins
$$(star_{K^\prime}(\hat{\sigma}),link_{K^\prime}(\hat{\sigma}))=
\partial \sigma^\prime * (D(\sigma,K),\partial D(\sigma,K))$$
\end{definition}
\textit{Proof of Poincar\'e duality.}
We will prove Poincar\'e duality following the proof for compact manifolds
in \cite{gh}. Let $A$ and $B$ be two cycles in a  PL
manifold $M$, intersecting transversely at $p$. Let $\iota_p(A \cdot B)$ be
the signed intersection of $A$ with $B$ at $p$ (which is either +1 or -1). 
Note that $D(\sigma)$ transversely intersects $\sigma$ at a point $p$.
We choose orientations of $\sigma$ and $D(\sigma)$ so that
$\iota_p(\sigma,D(\sigma))=+1$ for any simplex $\sigma$ in $K$. 

We now relate the boundary operator $\partial$ on the complex $\{\sigma_\alpha\}$ to
the coboundary operator $\delta$ on $\{D(\sigma_\alpha)\}$. Note first that if
$\sigma$ has vertices $\sigma^0,...,\sigma^k$, then the dual cell is given
by the $(k+1)$-fold intersection of the dual $n$-cells to the vertices.
So the cells appearing in the coboundary $\delta D(\sigma)$ will just be
the k-fold intersection of the dual cells of the faces of $\sigma$. We
have the relation
$$\delta(D(\sigma))=(-1)^{n-k+1}D(\partial \sigma)$$
(see \cite{gh} pp.54-55 for a proof).From this we see that the map 
$\sigma \to D(\sigma)$ induces an isomorphism between the complex
$(C^{uff}_*(K,\textbf{Z}),\partial)$ of uniformly finite chains of the original simplicial
decomposition of $M$ and the complex $(C_{uff}^*(K^\prime,\textbf{Z}),\delta)$ of
cochains in the dual cell decomposition. This proves Poincar\'e duality. 

\bigskip\noindent

We will also need the analogue of the Eilenberg-Zilber theorem for ordinary
uniformly homology and cohomology theories. This theorem is due to the author
and J.Block and also appears in \cite{ab}.

\begin{theorem}[Eilenberg-Zilber Theorem]
Let $C^{uff}_n(X;\textbf{Z})$ be the \break uniformly
finite cellular $n$-chains with coefficients in $\textbf{Z}$. Then
$$C_n^{uff}(X \times Y;\textbf{Z})=\bigoplus_{i+j=n} C_i^{uff}(X;C_j^{uff}(Y;\textbf{Z}))$$
where $C_i^{uff}(X;C_j^{uff}(Y;\textbf{Z}))$ means uniformly finite $i$-chains with values in the
uniformly finite $j$-chains. More explicitly, we take uniformly finite chains
with coefficients in the normed group $C^{uff}_j(Y;\textbf{Z})$, where the norm of a chain
is defined to be the supremum of its coefficients. The boundary map is given
by the Leibnitz rule.
\end{theorem}
\textit{Proof} Write $M=\{\sigma_\alpha^k\}_{\alpha,k}$ and $N=\{\sigma^{\prime
k}_\alpha\}_{\alpha,k}$. The products $\sigma_\alpha^k \times \sigma^{\prime l}
_\beta$ give a cell decomposition of the product $M \times N$ with the
boundary operator
$$\partial(\sigma^k_\alpha \times \sigma^{\prime l}_\beta)=\partial \sigma^k
_\alpha \times \sigma^{\prime l}_\beta + (-1)^k \sigma^k_\alpha \times
\partial \sigma^{\prime l}_\beta.$$

We compute the chains as cellular chains by the triangulations of $X \times Y$,
$X$ and $Y$. In general, we have:
$$(X \times Y)^{(n)}=\bigcup_j X^{(j)} \times Y^{(n-j)}$$
We note also that $l^\infty(S \times T)=l^\infty(S;l^\infty(T))$, where
$S$ and $T$ are discrete sets.

We claim that the homology of this chain complex is $H^{uff}_*(X \times Y)$.
We have a map
$$C_n^{uff}(M \times N;\textbf{Z}) \to \bigoplus_{i+j=n}C^{uff}_i(M;C_j^{uff}(
N;\textbf{Z}))$$
which we have shown to be an isomorphism. By the construction of the boundary,
this induces a map on homology. This clearly takes the boundary to the 
boundary, and is thus an isomorphism.

\section{The Whitehead Group}
We first recall some definitions regarding the $bg$ Whitehead group of a
simplicial complex of bounded geometry.

\begin{definition}
Let $X$ be a CW complex of bounded geometry. An expansion of
bounded geometry is a $bg$ CW complex $Y$ so that

i. $(Y,X)$ is a $bg$ CW pair.

ii. $Y=X \cup_f (\Sigma \times I^r) \cup_g (\Sigma \times I^{r+1})$ for
$bg$ $(r+i)$-cells $\Sigma \times I^{r+i}$, $i=0,1$ and attaching maps
$f,g$.

iii. There is a characteristic map $\psi_{r+1}:\Sigma \times I^{r+1} \to Y$
for the $bg$ $(r+1)$-cell so that $\psi_{r+1} \mid \Sigma \times I: \Sigma
\times I^r \to Y$ is characteristic for the $bg$ $r$-cell. If $Y$ is a $bg$
expansion of $X$, $Y$ is said to $bg$ collapse to $X$.
\end{definition}
\begin{definition}
Let $DR^{bg}$ be the collection of all pairs $(Y,X)$ so that
$X$ is a $bg$ strong deformation retract of $Y$. The collection of equivalence
classes of such a pairs under elementary $bg$ expansions and collapses rel $X$ 
is denoted $Wh^{bg}(X)$.
\end{definition}
We refer \cite{am} for the definitions of controlled bases, controlled modules
and $R\mathcal{P} G_1(X)$-modules.

\begin{definition}
Let $X$ be a metric space, $R$ a ring. A controlled basis is
a pro-object defined as follows. The basis is a pair $(S,\sigma)$ where $S$
is a set and $\sigma$ is a function from $S$ to the collection of open sets
in a control space $X$. A morphism from $(S, \sigma)$ to $(T,\rho)$ is given
by a map $\alpha:S \to T$ along with a natural transformation $\rho \alpha \to
C^n\sigma$, where $C^n$ is the operation on metric balls in the control space
which increases the radius by $n$.
\end{definition}
\begin{definition}
Let $(S, \sigma)$ be a controlled basis. The free
$R\mathcal{P} G_1(X)$-module $F(\sigma)$ with basis $(S,\sigma)$ is a functor
from $\mathcal{P} G_1(X)$ to the category of $R$-modules defined as follows:

i. For any $b$ an object in $\mathcal{P} G_1(X)$, $F(\sigma)(b)$ is the free 
$R$-module on $\{(\beta,s) \mid s \in S,\mbox{ } \beta \in \mathcal{P} G_1(X)(\sigma(s),b)\}$, where $\mathcal{P} G_1(X)(\sigma(s),b)$ is the set of morphisms
from $\sigma(s)$ to $b$ in $\mathcal{P} G_1(X)$.

ii. For any $\gamma \in \mathcal{P} G_1(X)(b,c)$, $\gamma_*:F(\sigma)(b) \to
F(\sigma)(c)$ has $\gamma_*(\beta,s)=(\gamma \beta,s)$.
\end{definition}
\begin{definition}
The category of controlled free $bg$ $\textbf{Z}\mathcal{P} G_1(X)$
modules is defined to be the category of controlled modules so that in any
ball of fixed radius the modules fall into a finite number of types. Morphisms
are defined to be morphisms of the modules so that if the control space is
partitioned into neighborhoods of a fixed radius, the restrictions fall into
a finite number of equivalence classes. By abuse of notation, we will denote
this category by $\textbf{Z}\mathcal{P} G_1(X)^{bg}$.
\end{definition}
\begin{definition}
We define $K_1(\textbf{Z}\mathcal{P} G_1(X)^{bg})$ to be the
abelian group generated by $[F,\alpha]$ where $F$ is a controlled free $bg$
module and $\alpha$ is an automorphism of $F$ so that

i. $[F,\alpha]=[F^\prime,\alpha^\prime]$ if there is an isomorphism
$\phi:F \to F^\prime$ so that $\phi\alpha=\alpha^\prime \phi$.

ii. $[F \oplus F^\prime, \alpha \oplus \alpha^\prime]=[F,\alpha]+[F,\alpha^\prime]$

iii. $[F,\alpha\beta]=[F,\alpha]+[F,\beta]$.
\end{definition}
\begin{definition}
$Wh^{bg}(\mathcal{P} G_1(X))$ is defined to be the quotient of
$K_1(\mathcal{P} G_1(X)^{bg})$ defined by taking the quotient by the subgroup of
elements of the form
$$[F(\sigma),u_{F(\sigma)}]$$
and
$$[F(\sigma),F(\alpha,\nu)]$$
where $(S,\sigma)$ is any $bg$ basis over $\mathcal{P} G_1(X)$, $u_{F(\sigma)}$
is multiplication by a unit, and $F(\alpha,\nu)$ is an automorphism of bases.
\end{definition}
\begin{definition}
Let $\mathcal{A}$ be a small additive category. The idempotent
completion $\hat{\mathcal{A}}$ is the category with objects morphisms $p:A \to A$
of $\mathcal{A}$ so that $p^2=p$ and morphisms given by morphisms $\phi:A_1 \to A_2$
so that $\phi=p_1\phi p_2$, where $p_i:A_i \to A_i$ are the source and
target.
\end{definition}
\begin{definition}
$K_0^{bg}(X)$ the $bg$ projective class group of $X$, where
$X$ is a $bg$ simplicial complex, is defined to be $K_0$ of the idempotent
completion of the category $\textbf{Z}\mathcal{P} G_1(X)^{bg}$. There is a
homomorphism
$$rank:K_0^{bg}(X) \to H_0^{uff}(X;\textbf{Z})$$
given as follows. Let $m \in K_0^{bg}(X)$ be a given element. Then one can
find a representative for $m$ which has basis elements only at each vertex
of the simplicial complex $X$. Thus to a given vertex one can naturally
associate a free module. Define rank($m$) to be the uniformly finite 0-chain
rank($m$)=$\sum r_x x$, where $x$ is the rank of the free module constructed
above. This clearly gives a map $K_0^{bg}(X) \to C_0^{uff}(X;\textbf{Z})$,
and we observe that by taking infinite process tricks into account which
cancel the ranks, we can pass to homology. The kernel of this map is the
reduced $bg$ projective class group of $X$ and is denoted $\tilde{K}^{bg}_0(X)$.
It measures the obstruction to a $bg$ projective module being free.
\end{definition}
\begin{definition}Let $X$ be a $bg$ simplicial complex. Define the group 
$K_{-i}^{bg}(X)$, for $i>0$, to be $K_{-i}(\textbf{Z}\mathcal{P}G_1(X)^{bg})$.
\end{definition}
The following theorem is proven in \cite{attie}:
\begin{theorem} There is an isomorphism between the geometric Whitehead \break group,
and the algebraic Whitehead group
$$Wh^{bg}(M) \simeq Wh^{bg}(\mathcal{P}G_1(M)).$$ 
\end{theorem}
Before stating the $bg$ s-cobordism theorem, we remark on the necessary notions
of handlebody theory. It is a classical fact that one can construct a 
handlebody from a triangulation of a PL manifold, so that Theorem 2.1
yields a handlebody
decomposition for a $bg$ PL manifold. To prove the $bg$ s-cobordism theorem,
which will relate $bg$ h-cobordisms to the $bg$ controlled Whitehead group, we
need to introduce a $bg$ controlled handlebody decomposition. This is identical
to the controlled handlebody decomposition of \cite{am}, except that one uses
the triangulation introduced in section 2 in place of the usual triangulation.
Furthermore, all of the handle constructions and operations work identically
as in \cite{am} except that they must be carried out in a uniform manner.
$bg$ transversality must also be used in place of ordinary transversality 
everywhere. We omit the details.
\begin{theorem}[s-cobordism theorem] Let $M$ be a PL manifold of bounded 
geometry controlled by a $bg$ map over a uniformly contractible space $Z$.
Then $Wh^{bg}(M)$ is in one-to-one correspondence with set of h-cobordisms
over $M$, whenever $M$ is a PL manifold of bounded geometry of dimension
$\ge 5$. In particular, if $W$ is a $bg$ h-cobordism which is $bg$ simple 
homotopy equivalent to one end. Then it is $bg$ PL homeomorphic to $M \times I$.\end{theorem}
\begin{theorem}[Smooth s-cobordism theorem] A smooth $bg$ s-cobordism  \break of 
dimension $\ge 6$ is $bg$ diffeomorphic to a product.
\end{theorem}

\section{Surgery Theory}
The first four definitions are taken from \cite{attie}. They are all based on
the controlled surgery theory of Ferry-Pedersen \cite{fp}. See also \cite{hp}.
\begin{definition}
A Poincar\'e duality space $Y$ of bounded geometry over a
simplicial complex $X$ of bounded geometry is defined to be a simplicial
complex of bounded geometry over $X$ so that there is a fundamental class
$[Y]$ in the top dimensional locally finite homology
$$[Y] \in H_n^{lf}(Y;\textbf{Z})$$
of $Y$ so that taking the
cap product $[Y]\cap -:D^\#(Y) \to D_{n-\#}(Y)$ induces a homotopy equivalence
of controlled chain complexes. $Y$ is a simple $bg$ Poincar\'e duality
space if the torsion of $[Y] \cap -$ is trivial in $Wh^{bg}(Y)$.
\end{definition}
\begin{definition} Let $p:(Y,\partial Y) \to X$ be a proper bounded CW
pair so that $X$ has bounded fundamental group $\pi$. The pair $(Y, \partial Y)$
is an n-dimensional $bg$ Poincar\'e duality pair if $\partial Y$ is an
$(n-1)$-dimensional $bg$ Poincar\'e complex with orientation double covering
the pullback of the orientation double covering on $X$ and if there is an
element $[Y] \in H^{lf}_n(Y,\partial Y; \textbf{Z})$ such that
$[Y]\cap -:D^\#(X) \to D_{n-\#}(X,\partial X)$ is a homotopy equivalence
of chain complexes. $(Y,\partial Y)$ is a simple $bg$ Poincar\'e duality
space, if the torsion of $[Y] \cap -$ is trivial in $Wh^{bg}(Y)$.
\end{definition}
\begin{definition}
A $bg$ spherical fiber space $E$ over a $bg$ simplicial
complex $X$ is a bundle $p:E \to X$ with fiber a sphere $S^k$ so that $E$
is $bg$ simplicial complex and $p$ is a $bg$ map.
\end{definition}
The following definition is from \cite{fp}:
\begin{definition} Let $(X,p)$ be a bounded Poincar\'e duality space.
Construct a proper embedding $X \subset \textbf{R}^n$, $n-dimX \ge 3$. Let
$W$ be a regular neighborhood of $X$ and $r: W \to X$ a retraction. $W \to M$
has a bounded fundamental group, and we can triangulate sufficiently finely
to get a bounded CW structure on $W$. Then the controlled Spivak normal
fibration is the fibration $\partial W \to X$. Let $F$ be the homotopy fiber
of this fibration.
\end{definition}
\newtheorem{lemma}{Lemma}[section]
\begin{lemma}[\cite{fp}] The fibre $F$ is homotopy equivalent to sphere of
dimension $n-dimX-1$.
\end{lemma}
\begin{definition}
Let $X$ be a $bg$ Poincar\'e duality complex. The $bg$ Spivak normal fibration
is defined by giving the controlled Spivak fibration of $X$ a $bg$ structure by
observing that the projection map can be taken to be boundedly controlled
and hence is $bg$.
\end{definition}
\begin{definition} Let $X^n$ be a $bg$ Poincar\'e duality space over a
$bg$ simplicial complex $M$ and let $\nu$ be a $bg$ $PL$ block bundle over $X$. 
A bounded geometry surgery problem is a triple $(W^n,\phi,F)$ where
$\phi:W \to X$ is a $bg$ map from an $n$-manifold $W$ to $X$ such that
$\phi_*([W])=[X]$ and $F$ is a stable trivialization of $\tau_W \oplus
\phi^*\nu$. Two problems $(W, \phi, F)$ and $(\overline{W},\overline{\phi},
\overline{F})$ are equivalent if there exist an $(n+1)$-dimensional manifold
$P$ with $\partial P=W \coprod \overline{W}$, a $bg$ map $\Phi:P \to X$
extending $\phi$ and $\overline{\phi}$, and a stable trivialization
of $\tau_P \oplus \Phi^*\nu$ extending $F$ and $\overline{F}$.
\end{definition}
\begin{definition}
The $bg$ surgery group of a simplicial complex of bounded
geometry is defined as follows. An unrestricted object consists of:

1. A $bg$ Poincar\'e pair $(Y,X)$ over a $bg$ complex $M$.

2. A $bg$ map $\phi:(W,\partial W) \to (Y,X)$ of pairs of degree 1, where
$W$ is a PL manifold of bounded geometry and $\phi \mid:\partial W \to X$
is a $bg$ simple homotopy equivalence.

3. A $bg$ stable framing $F$ of $\tau_W \oplus \phi^*(\tau)$, where $\tau$
is the $bg$ Spivak normal fibration of $(Y,X)$.

4. A map $\omega:Y \to K$, where $K$ is $bg$ complex so that the pullback of
the double cover of $K$ to $Y$ is orientation preserving.

The surgery group is then defined to be the cobordism group of such unrestricted
objects. It is denoted $L^{bg}_n(K).$
\end{definition}
Unrestricted objects, however, cannot be used for surgery. To take care of this
one introduces restricted objects, for which the map $\omega$ induces an
isomorphism of the fundamental groups.

The next definition is based on \cite{am}.
\begin{definition} Let $T$ be a discrete set. A $bg$ $r$-handle is a $bg$ pair
$(T \times (D^r,S^{r-1}) \times D^{n-r})$ with a $bg$ control map over a 
$bg$ complex $Z$. A uniform surgery is an exchange of $T_1 \times S^r \times
D^{n-r}$ by $T_1 \times D^{r+1} \times S^{n-r-1}$ where $T_1$ is a discrete 
set, and the regluing PL homeomorphisms 
along each $S^r \times S^{n-r-1}$ in $T_1 \times S^r
\times S^{n-r-1}$ lie in a finite set.
\end{definition}
We recall the following result from \cite{attie}:
\begin{lemma}
Let $(X^n,\partial X)$ be a $bg$ Poincar\'e duality space over
$M$, $n \geq 6$. Consider a $bg$ surgery problem $\phi:(W,\partial W) \to
(X, \partial X)$. Then $\phi$ is equivalent to a $bg$ surgery problem
$\overline{\phi}:(\overline{W}, \partial \overline{W}) \to (X, \partial X)$
so that $\phi$ is $bg$ $[\frac{n}{2}]$ connected over $M$ and is $[\frac{n-1}{2}]$
connected when restricted to the boundary.
\end{lemma}
The following theorem is proven analogously to the corresponding $\pi-\pi$
theorem in \cite{fp}.

\begin{theorem}
Let $\phi$ be a $bg$ map of pairs $\phi:(X,Y) \to (N,M)$. Suppose
the inclusion $Y \subset X$ induces an isomorphism on the controlled
fundamental groups. Then one can uniformly surger $\phi$ to obtain a homotopy
equivalence.
\end{theorem}
\textit{Proof.}
By Lemma 5.4 of \cite{fp} we may do surgery
up to the middle dimension. This means that cancelling cells in the controlled
algebraic mapping cone of the corresponding controlled chain complexes
$$D_{\#}(W^\prime, \partial W^\prime ; W^\prime, \partial W^\prime)$$
yields a complex which is 0 through dimension $[\frac{n}{2}]$. Here
$$W^\prime \overset{\phi^\prime}\longrightarrow X$$
is the surgery problem obtained so that $\phi^\prime$ is an inclusion which
is the identity below the middle dimension. The $k+1$ dimensional generators
are represented by $k$-dimensional discs $D$ in $\partial X$ whose boundaries
lie in $\partial W^\prime$. Note that there is a parallel copy of $D$ in a
collar neighborhood of $\partial W^\prime$ which is contained in $W^\prime$.
Now use cell trading to change $D_{\#}$ to
$$0 \to D^\prime_{k+3} \overset{\partial} \longrightarrow D^\prime_{k+2}
\overset{\partial}\longrightarrow D_{k+1} \to 0$$
together with a homotopy $s$ so that $s\partial+\partial s=1$ except at degree
$k+1$. Corresponding to each generator of $D_{k+2}$ we introduce a pair of
cancelling $(k-1)$- and $k$-handles and excise the interior of the $(k-1)$-handle
from $(W, \partial W)$. The modified chain complex is:
$$0 \to D_{k+3} \to D_{k+2} \to D_{k+1} \oplus D_{k+2} \to 0$$
All generators of $D_{k+1}\oplus D_{k+2}$ are represented by discs. We may
represent any linear combination of these discs by an embedded disc, and these
embedded discs may be assumed disjoint by a piping argument which uses the
$\pi-\pi$ condition in the hypothesis. We do surgery on the following
elements: For each generator $x$ of $D_{k+1}$, we do surgery on
$(x-s\partial x, sx)$ and for each generator $y$ of $D_{k+3}$, we do surgery
on $(0,\partial y)$. This results in a contractible chain complex and completes
the even dimensional case. The odd dimensional case follows by crossing
with $S^1$ and splitting back. This can be done by a codimension 1 splitting
technique which is simpler than Theorem 5.5 of \cite{attie} and therefore
will not be worked out in detail.

The following is a direct corollary of the above, as in \cite{wall}.

\begin{theorem}
Let $\phi:(W, \partial W) \to (Y,X)$ be an unrestricted object
which is $bg$ 1-connected and a $bg$ homotopy equivalence on the boundary.
Then there is a $bg$ normal cobordism rel $\partial W$ to a $bg$ homotopy
equivalence if and only if the equivalence class of $\phi$ in $L^{bg}(K)$
vanishes.
\end{theorem}

\begin{definition}
Define $NI^{bg,PL}(X)$ to be cobordism group of triples
$(M,\phi,F)$, where $M$ is a PL manifold of bounded geometry, $\phi$ a degree
one $bg$ normal map to $X$ and $F$ a stable $bg$ trivialization of
$\tau_M \oplus \phi^*\nu$, where $\nu$ is the $bg$ normal bundle of $X$.
\end{definition}
\begin{definition}
Define the simple $bg$ PL structure set $\mathcal{S}^{bg,s}_{PL}(X)$, where $X$ is
a PL manifold of bounded geometry to be the set of $bg$ simple homotopy
equivalences $\phi:N \to X$ modulo the equivalence relation $\phi \sim \phi^\prime:
N^\prime \to X$, if there is a PL quasi-isometry $h:N \to N^\prime$ so that
$\phi^\prime \circ h=\phi$.

Define the simple $bg$ TOP structure set $\mathcal{S}^{bg,s}_{TOP}(X)$, where
$X$ is a PL manifold of bounded geometry to be the set of $bg$ simple
homotopy equivalences $\phi:N \to X$ modulo the equivalence relation
$\phi \sim \phi^\prime: N^\prime \to X$, if there is a continuous quasi-isometry
$h:N \to N^\prime$ so that $\phi^\prime \circ h=\phi$.
\end{definition}
\begin{proposition}
There is an exact sequence
$$\mathcal{S}^{bg,s}_{PL}(X) \to NI^{bg,PL}(X) \to L^{bg}_n(X)$$
\end{proposition}

We next introduce an algebraic version of this theory, following Ranicki. We
recall that Ranicki has introduced the L-theory of an additive category. We
refer to Ranicki's book \cite{ranicki1} for the relevant definitions. The
following treatment of algebraic surgery is based directly and heavily on
\cite{ranicki1}.

The following theorem relates the surgery groups defined above to the L-theory
of an additive category. 
\begin{definition}
Let $K$ be a simplicial complex of bounded geometry. Let
$\mathcal{A}$ be an additive category, $\pi$ a group. Then the category
$\mathcal{C}^{bg}_K(\mathcal{A}[\pi])$ is defined to be the one whose objects are formal
direct sums
$$M=\sum_{x \in K} M(x)$$
of objects $M(x)$ in $\mathcal{A}[\pi]$, which fall into a fixed finite number of
types inside of each ball of fixed radius in $K$. Here $\mathcal{A}[\pi]$ is the
category with the one object $M[\pi]$ for each object $M$ in  $\mathcal{A}$, and with morphisms linear
combinations of morphisms $f_g:M \to N$ in $\mathcal{A}$ of the form
$$f=\sum_{g \in \pi} n_g f_g: M[\pi] \to N[\pi],$$
with $\{g \in \pi \mid f_g \neq 0 \}$ finite.
We use the notation $\mathcal{C}^{bg}_i(\mathcal{A}[\pi])$ for 
$\mathcal{C}^{bg}_{\textbf{R}^i}(\mathcal{A}[\pi])$. For any commutative ring
$R$ there is an identification
$$\mathcal{A}^h(R)[\pi]=\mathcal{A}^h(R[\pi])$$
with $\mathcal{A}^h(R)$ the additive category of based finitely generated
free $R$-modules. Write the category $\mathcal{C}_X^{bg}(\mathcal{A}^h(R[\pi]))$
as $\mathcal{C}_X^{bg}(R[\pi])$. 
\end{definition}

\begin{theorem} Let $X$ be a manifold of bounded geometry with bounded
fundamental group $\pi_1(X)$. Then
$L_*^{bg}(X)=L_*(\mathcal{C}^{bg}_X(\textbf{Z}\pi_1(X)))$. Moreover every
element of $L_*(\mathcal{C}^{bg}_X(\textbf{Z}\pi_1(X)))$ is realized as the
obstruction on a surgery problem with target $N \times I$ and homotopy
equivalence on the boundary for an arbitrary $n-1$ dimensional $bg$
PL manifold $N$ with bounded fundamental group $\pi_1(X)$. 
\end{theorem}
\textit{Proof.} This is proven as in \cite{fp}. In the even dimensional case
we first obtain a highly connected surgery problem. We obtain a chain
complex homotopy equivalent to $K_{\#}(M)$ which is concentrated in dimensions
$k+2$, $k+1$ and $k$, and a contracting homotopy $s$ which is obtained from
Poincar\'e duality. Introducing cancelling $k+1$ and $k+2$ handles, we may
shorten this chain complex to a 2-term chain complex
$$0 \to K^\prime_{k+1} \to K_k^\prime \to 0$$
We can then do further surgery to get a chain complex concentrated in one
degree. Denote the remaining module by $A$. Poincar\'e duality produces an
isomorphism $\phi:A \to A^*$ which determines the intersections of different
generators, i.e. $\phi(e_i)(e_j)$ determines the intersections of $e_i$ and
$e_j$ when $e_i$ and $e_j$ are different. Now total order the basis and
define a map $\nu:A \to A^*$ so that $\nu(e_i)(e_j)$ is 0 when $i>j$ and
the intersection counted with sign in $\textbf{Z}\pi_1(X)$ when $i \leq j$.

By symmetrization $\nu+\epsilon\nu^*=\phi$, hence an isomorphism. This
represents the surgery obstruction. Using -1 and 0-connectedness, if the
surgery obstruction is zero, we may find a Lagrangian so that doing uniform surgery
on this Lagrangian produces a $bg$ homotopy equivalence.

In odd dimensions we do surgery below the middle dimension, and proceeding
as above we may obtain a length 2 chain complex
$$0 \to K_{k+1} \to K_{k} \to 0$$
Now do surgeries on embedded $S^k \times D^{k+1}$'s so that, denoting the
trace of the surgery by $W$, the chain complexes $K_{\#}(W,M)$,$K_{\#}(W)$
and $K_{\#}(W,M^\prime)$ are homotopy equivalent to chain complexes which
are zero except in dimension $k+1$. One way to do this could be to do
surgeries to all the generators of $K_k$. Denote the resulting manifold
by $M^\prime$. The surgery obstruction is now defined to be the following
formation
$$(K_{k+1}(W,M) \oplus K_{k+1}(W,M^\prime), K_{k+1}(W,M), K_{k+1}(W))$$
where the first Lagrangian is the inclusion on the first factor, and the
second Lagrangian is induced by the pair of inclusions. Poincar\'e duality
shows that these are indeed Lagrangians. This is a well-defined element in
the odd L-group.

Realization also follows as in \cite{fp}. 
 
Let $\mathcal{B}(\mathcal{A})$ be the additive category of finite chain complexes in
$\mathcal{A}$ and chain maps.

\begin{definition}
A subcategory $\mathcal{C} \subseteq \mathcal{B}(\mathcal{A})$ is closed if
it is a full additive subcategory so that the algebraic mapping cone $C(f)$
of any chain map $f:C \to D$ in $\mathcal{C}$ is an object of $\mathcal{C}$. A chain
complex $C$ in $\mathcal{A}$ is $\mathcal{C}$-contractible if it belongs to $\mathcal{C}$.
A chain map $f:C \to D$ in $\mathcal{A}$ is a $\mathcal{C}$-equivalence if $C(f)$ is
$\mathcal{C}$-contractible.

An $n$-dimensional quadratic complex $(C,\psi)$ in $\mathcal{A}$ is $\mathcal{C}$-contractible
if $C$ and $C^{n-*}$ are $\mathcal{C}$-contractible. An $n$-dimensional quadratic
complex $(C,\psi)$ in $\mathcal{A}$ is $\mathcal{C}$-Poincar\'e if the chain complex
$$\partial C = S^{-1}C((1+T)\psi_0:C^{n-*} \to C)$$
is $\mathcal{C}$-contractible.
\end{definition}
\begin{definition}
Let $\Lambda=(\mathcal{A}, \mathcal{B}, \mathcal{C})$ be a triple of additive
categories, where $\mathcal{A}$ has a chain duality $T: \mathcal{A} \to \mathcal{B}(\mathcal{A})$
and a pair $(\mathcal{B}, \mathcal{C} \subseteq \mathcal{B})$ of closed subcategories of
$\mathcal{B}(\mathcal{A})$ so that for any object $B$ of $\mathcal{B}$

i. The algebraic mapping cone $C(1:B \to B)$ is an object of $\mathcal{C}$

ii. The chain equivalence $e(B):T^2(B) \to B$ is a $\mathcal{C}$-equivalence.

Then $\Lambda$ is said to be an algebraic bordism category.
\end{definition}
\begin{definition} For any additive category with chain duality
$\mathcal{A}$ there is defined an algebraic bordism category
$$\Lambda(\mathcal{A})=(\mathcal{A}, \mathcal{B}(\mathcal{A}),
\mathcal{C}(\mathcal{A}))$$
with $\mathcal{B}(\mathcal{A})$ the category of finite chain complexes in
$\mathcal{A}$, and $\mathcal{C}(\mathcal{A}) \subseteq \mathcal{B}(
\mathcal{A})$ the subcategory of contractible complexes.
\end{definition}
\begin{definition}
Let $K$ be a simplicial complex. Let $\Lambda=(\mathcal{A},
\mathcal{B}, \mathcal{C})$ be an algebraic bordism category. An $n$-dimensional
quadratic complex $(C,\psi)$ in $\Lambda$ is an $n$-dimensional quadratic
complex in $\mathcal{A}$ which is $\mathcal{B}$-contractible and $\mathcal{C}$-Poincar\'e.
The quadratic $L$-group $L_n(\Lambda)$ is the cobordism group of $n$-dimensional
quadratic complexes in $\Lambda$.
\end{definition}
\begin{definition}
Let $K$ be a simplicial complex. An object $M$ in an additive
category $\mathcal{A}$ is said to be $K$-based if it is expressed as a direct sum
$$M=\sum_{\sigma \in K} M(\sigma)$$
of objects $M(\sigma)$ in $\mathcal{A}$ so that $\{\sigma \in K: M(\sigma) \neq 0\}$
is finite in ball of fixed radius in $K$. A morphism $f:K \to N$ of $K$-based
objects is a collection of morphisms in $\mathcal{A}$
$$f=\{f(\tau,\sigma):M(\sigma) \to N(\tau): \sigma, \tau \in K\}$$
\end{definition}
\begin{definition}
Let $K$ be a simplicial complex of bounded geometry. Let \break
$\mathcal{A}^{uf}(K)$ be the additive category of $K$-based objects in $\mathcal{A}$
which fall into a finite number of types in each ball of fixed radius in $K$.
Define a uniformly finite assembly map
$$\mathcal{A}^{uf}(K) \to \mathcal{C}^{bg}_K(\mathcal{A}[\pi])$$
by associating to $M$ the $bg$ controlled module $\hat{M}$ which has the  value
$M(\sigma)$ at the barycenter of $\sigma$ and is 0 everywhere else.
\end{definition}
\begin{definition}
Let $\mathcal{A}(R)^{uf}_*(K)$ be the additive category of
$K$-based objects in $\mathcal{A}(R)$ which fall into a fixed number of types in
a ball of fixed radius, with morphisms $f:M \to N$ such that $f(\tau,\sigma)=0:
M(\sigma) \to N(\tau)$ unless $\tau \geq \sigma$ so that $f(M(\sigma))
\subseteq \sum_{\tau \geq \sigma} N(\tau)$.
\end{definition}
\begin{definition}
 Define three algebraic bordism categories:

i. local, uniformly finite, finitely generated free $(R,K)$-modules
$$\Lambda(R)^{uf}_*(K)=(\mathcal{A}^{uf}(R,K),\mathcal{B}^{uf}(R,K),\mathcal{C}(R)^{uf}_*(K))$$
where $\mathcal{B}^{uf}(R,K)$ is the category of finite chain complexes of
f.g. free uniformly finite $(R,K)$-modules. An object in $\mathcal{C}^{uf}(R)_*(K)$
is finite f.g. free uniformly finite $(R,K)$-module chain complex $C$ such that each
$[C][\sigma]$ ($\sigma \in K$) is a contractible f.g. free $R$-module chain
complex. 

ii. global, uniformly finite, finitely generated free $(R,K)$-modules
$$\Lambda^{uf}(R,K)=(\mathcal{A}^{uf}(R,K),\mathcal{B}^{uf}(R,K),\mathcal{C}^{uf}(R,K))$$
with
$$\mathcal{C}^{uf}(R,K) \subseteq \mathcal{B}^{uf}(R,K)$$ the subcategory of finite
finitely generated free uniformly finite $(R,K)$-module chain complexes
$C$ which assemble to contractible finitely generated free \break
$\mathcal{C}^{bg}_K(R[\pi])$-module chain complexes.

iii. $\mathcal{C}^{bg}_K(\mathcal{A}[\pi])$, along with the categories of
chain complexes and contractible chain complexes over it. 
\end{definition}
\begin{definition}
The quadratic $bg$ structure groups of $(R,K)$ are the
cobordism groups
$$\mathcal{S}^{bg}_n(R,K)=L_{n-1}(\mathcal{A}^{uf}(R,K), \mathcal{C}^{uf}(R,K),\mathcal{C}^{uf}(R)_*(K)).$$
\end{definition}
\begin{definition}
Let $M$ be a manifold of bounded geometry. We define an
object which we will refer to as uniformly finite homology with coefficients
in the L-spectrum,
$$H_n^{uff}(M;\textbf{L})=L_n(\Lambda(\textbf{Z})^{uf}_*(M))$$
For the moment this is a purely formal definition. At the end of this paper
we will show that for $M=N \times \textbf{R}^n$, $N$ a compact manifold, this
object is in fact the uniformly finite homology with coefficients in the
L-groups.
\end{definition}
\begin{theorem} A normal map of $n$-dimensional PL manifolds of bounded
geometry $(f,b):N \to M$ determines an element, the normal invariant
$$[f,b] \in H_n^{uff}(M;\textbf{L})$$
which assembles to the surgery obstruction. 
\end{theorem}
\textit{Proof.} Let $X$ be the polyhedron of an $n$-dimensional geometric $bg$ Poincar\'e
complex with a homotopy equivalence $g:M \to X$ so that both $g$ and
$gf:N \to X$ are $bg$ transverse to the dual cell decomposition
$\{D(\tau,X):\tau \in X\}$ of $X$. The restrictions of $f$ define a $uf$
cycle of degree 1 $bg$ normal maps of $(n-\mid \tau \mid)$-dimensional
manifolds with boundary
$$\{(f(\tau),b(\tau))\}:\{N(\tau)\} \to \{M(\tau)\}$$
with
$$M(\tau)=g^{-1}D(\tau,X), N(\tau)=(gf)^{-1}D(\tau,X), \tau \in X$$
so that $M(\tau)=\{pt.\}$ for $n$-simplices $\tau \in X^{(n)}$. The controlled
kernel $uf$ cycle
$$\{(C(f(\tau)^!),\psi(b(\tau))): \tau \in X\}$$
of $(n-\mid \tau \mid)$-dimensional quadratic Poincar\'e pairs $\mathcal{A}(\textbf{Z})$
is a 1-connective, $n$-dimensional quadratic Poincar\'e complex in
$\mathcal{A}(\textbf{Z})^{uf}_*(M)$ allowing the definition
$$[f,b]=\{(C(f(\tau)^!),\psi(b(\tau)))\} \in L_n(\Lambda(\textbf{Z}^{uf}_*(M)))$$

\bigskip\noindent

Finally, we need the uniformly finite version of the algebraic $\pi-\pi$
theorem.

\begin{theorem}
The global $uf$ assembly maps define isomorphisms
$$L_n(\mathcal{A}^{uf}(R,K)) \simeq L_n(\mathcal{C}_K^{bg}(\mathcal{A}[\pi_1(K)]))$$
\end{theorem}
\textit{Proof.} The proof is virtually the same as the one given by Ranicki in
\cite{ranicki1} in the
compact case, so we will restrict ourselves to those points where the two
proofs differ. In the compact case, the Hurewicz theorem is used to compare
the $(R,K)$-module with its image under the assembly. In the $bg$ case, this
Hurewicz theorem is replaced by the controlled Hurewicz theorem of Anderson
and Munkholm. The rest is verbatim the same as the one in Ranicki except that
modules over the entire space are replaced by controlled modules.

\bigskip\noindent
\begin{theorem}[Algebraic Surgery Exact Sequence]Let $K$ be a $bg$ simplicial
complex. There is an exact sequence
$$L_{n+1}(\mathcal{C}^{bg}_K(\textbf{Z}\pi_1(K)) \to \mathcal{S}^{bg}_{n+1}(\textbf{Z},K) \to H_n^{uf}(K;\textbf{L}) \to L_n(
\mathcal{C}^{bg}_K(\textbf{Z}\pi_1(K)))$$
\end{theorem}
\textit{Proof.} This is simply the exact sequence of Ranicki \cite{ranicki1},
Proposition 3.9.

\begin{theorem} The algebraic $bg$ structure group $\mathcal{S}^{bg}_n(\textbf{Z},K)$
of a uniform triangulation $K$ of a $bg$ PL manifold $M$, is isomorphic
to the topological structure set $\mathcal{S}_{TOP}^{bg,s}(M)$, with
group structure given by characteristic variety addition.
\end{theorem}
\textit{Proof.} As in \cite{ranicki1} p.198,
if $f:N \to M$ is a homotopy equivalence the quadratic complex
giving the normal invariant of $f$ is globally contractible, allowing the
definition
$$s(f)=\{C(f(\tau)^!,\psi(b(\tau)))\} \in \mathcal{S}_{n+1}(M)$$
This proves the result.

\bigskip\noindent

We wish to express the PL $bg$ surgery exact sequence as a version of Ranicki's
exact sequence. In order to do this we need to incorporate the Casson-Sullivan
obstruction to the uniqueness of triangulations. The triangulation of the
polyhedron of a $bg$ PL manifold is not unique. In other words, there are
$bg$ homeomorphic $bg$ PL manifolds, which are not $bg$ PL homeomorphic. Any
$bg$ homeomorphism of polyhedra can be approximated by a $bg$ PL map, however
this map does not necessarily have a $bg$ PL inverse. This leads to the
non-uniqueness of the $bg$ PL triangulation, which means that there are
non-equivalent (via a $bg$ PL homeomorphism) triangulations (or $bg$ PL types)
in the $bg$ homeomorphism type of a $bg$ PL manifold. One can show by
obstruction theory that, if one disregards the $n-3$ skeleton, the triangulation
on a $bg$ PL manifold is unique up to $bg$ PL homeomorphism, within that
$bg$ homeomorphism type. However, there is an ambiguity when one gets to the
$n-3$ skeleton caused by the fact that the signature of a topological 4-manifold
is divisible by 8, whereas the signature of a PL or smooth manifold is
divisible by 16. (See pp.182-183 of \cite{ranicki1}.) The ambiguity is resolved
by the Casson-Sullivan obstruction, which is a local obstruction defined
following \cite{ranicki1}. Let $f:M \to M^\prime$ be a $bg$ homeomorphism,
$(W^{n+1}; M, M^\prime)$ a $bg$ PL cobordism with $F \mid_M=id$, $F\mid_{M^\prime}=
f$
$$(F,B):(W^{n+1};M,M^\prime) \to M \times ([0,1];\{0\},\{1\})$$
and let $\sigma_*(F,B)=(C,\psi)$ be the $bg$ surgery obstruction of $(F,B)$.
Then we define an invariant
$$\kappa(f)=\sum_{\sigma \in M^{(n-3)}} (\mbox{signature}(C(\sigma),\psi(\sigma))/8)\sigma$$
which is an element of $H^3(M;\textbf{Z}_2)$. That this classifies triangulations
follows exactly as in the compact case. The Casson-Sullivan invariant
clearly defines a map $L^{bg}_{n+1}(M) \oplus H_{n-3}^{lf}(M;\textbf{Z}_2)
\to \mathcal{S}^{bg,s}_{PL}(M)$ by assigning the homotopy equivalence
$f:M^\prime \to M$ to the image of the surgery problem $(F,B):W \to M \times [0,1]$.

We also need to define the Kirby-Siebenmann invariant. This lies in \break
$H^{lf}_{n-4}(M;\textbf{Z}_2)$. If we let $(g,c):N \to M$ be an EPL degree
one normal map, where $N$ is a topological manifold which is a $bg$ Poincar\'e
duality complex controlled over $M$, we can define the surgery obstruction of
$(g,c)$, $\sigma_*(g,c)$ in terms of the cellular chains of $N$ as in \cite{fp}
and note that this defines a $bg$ surgery obstruction $\sigma_*(g,c) \in
L_n^{bg}(M)$. If we take the image of this obstruction in $H_{n-4}^{lf}(M;\textbf{Z}_2)$,
as with the Casson-Sullivan invariant above, we obtain the Kirby-Siebenmann
obstruction. This can be written
$$g_*\kappa(N)=\sum_{\sigma \in M^{(n-4)}}(\mbox{signature}(C(\sigma),\psi(\sigma))/8)\sigma$$
where $\sigma_*(g,c)=(C,\psi)$. 

We obtain from these theorems the following algebraic exact sequence which
was suggested to the author by A.Ranicki:

\begin{theorem}
Let $M$ be a manifold of bounded geometry. Then there is an
exact sequence
$$\to L^{bg}_{n+1}(M) \oplus H^{lf}_{n-3}(M;\textbf{Z}_2) \to \mathcal{S}^{bg,s}_{PL}(M)
\to H_n^{uff}(M;\textbf{L})$$$$ \to L_n^{bg}(M) \oplus H_{n-4}^{lf}(M;\textbf{Z}_2)$$
\end{theorem}
\textit{Proof.}
We now prove exactness at $\mathcal{S}^{bg,s}_{PL}(M)$. A $bg$ homotopy equivalence
$f:M^\prime \to M$ has zero normal invariant in $H_n^{uff}(M; \textbf{L})$
if there is a $bg$ PL normal bordism $(F^\prime, B^\prime):(W;M^\prime, M^{''})
\to M \times ([0,1];\{0\},\{1\})$ to a $bg$ homeomorphism $f^{''}=F^\prime \mid
M^{''}:M^{''} \to M$. Hence $W$,$F^\prime$, $B^\prime$,$f^{''}$
are the same as $W$, $F$, $B$ and $f$ in the definition of the Casson-Sullivan
invariant above. Thus the kernel of the map $\mathcal{S}^{bg,s}_{PL}(M) \to
H_n^{uff}(M;\textbf{L})$ is equal to the image of the map defined above in
the definition of the Casson-Sullivan invariant. 
That is in this case, $(M,f)$ is in the image of $(\sigma_*(F^\prime,
B^\prime), \kappa(f^{''}))$. This proves exactness at $\mathcal{S}^{bg,s}_{PL}(M)$.

Next we prove exactness at $H_n^{uff}(M;\textbf{L})$. 
The kernel of the map 
$$H_n^{uff}(M;\textbf{L}) \to L_n(M)\oplus H^{lf}_{n-4}(M;
\textbf{Z}_2)$$
is equal to the set of degree one normal maps $f:M^\prime \to M$ with zero
surgery obstruction and zero Kirby-Siebenmann invariant. 
Suppose that $f:M^\prime \to M$ is a degree
one normal map with zero surgery obstruction and zero Kirby-Siebenmann
obstruction. Then by the definition, $f$ is a $bg$ PL structure on $M$. This
proves exactness at $H_n^{uff}(M;\textbf{L})$.

Finally, we prove exactness at 
$H_{n-3}^{lf}(M;\textbf{Z}_2) \oplus L^{bg}_{n+1}(M)$. The image of the map
$H_{n+1}^{uff}(M;\textbf{L}) \to L_{n+1}(M) \oplus H^{lf}_{n-3}(M;\textbf{Z}_2)$ 
is the assembly of the normal invariants of $(F,B):N \to M \times I$ to 
 surgery obstructions, which by the usual proof of the surgery exact
sequence of Wall \cite{wall} is the kernel of the map $L^{bg}_{n+1}(M) \oplus
H^{lf}_{n-3}(M; \textbf{Z}_2)\to \mathcal{S}^{bg,s}_{PL}(M)$.

\bigskip\noindent
Finally, we need to prove Siebenmann periodicity. This will follow immediately
from the above exact sequence. Note first that if $\sigma^*(f)$ is the
surgery obstruction of $f:M \to N$ then
$$\sigma^*(f)=\sigma^*(f \times CP^2)$$
This defines an isomorphism of the $bg$ L-groups $L^{bg}_m(M) \simeq L^{bg}_{m+4}(M)$.
Next we define the $bg$ ``resolution obstruction''. Let $(C \to D, (d\psi, \psi))$
be an $n$-dimensional locally Poincar\'e globally contractible quadratic pair
in $\mathcal{A}^{uf}(\textbf{Z},X)$ with $C$ 1-connective and $D$ 0-connective,
then the image of the algebraic complex
$$x=\sum_{r \in X^{(n)}}((D/C)(\tau),(\delta\psi/\psi)(\tau))\tau \in
H^{uff}_n(X;L_0(\textbf{Z}))=\textbf{Z}$$
in $\mathcal{S}^{bg,s}_{TOP}(X)$ is the ``resolution obstruction''. Note that this obstruction
is necessarily an element of \textbf{Z}. It is clear that we then have the
following standard periodicity result (cf. \cite{siebenmann2},\cite{ks},\cite{cappellw}, \cite{nicas},
\cite{weinbergeryan},\cite{sullivan})

\begin{theorem}
 The topological structure set of a manifold $M$ of bounded geometry is almost
4-fold periodic:
$$\mathcal{S}^{bg,s}_{TOP}(M) \simeq \mathcal{S}^{bg,s}_{TOP}(M \times D^4, \partial) \oplus \textbf{Z}$$
\end{theorem}
\section{Whitehead Group of $M \times \textbf{R}^n$}
The technique we use to calculate the Whitehead group is the same as used in
\cite{attie}, but more elaborate due to the fact that there are non-compact
directions involved. We will introduce a filtration of $Wh^{bg}(M \times
\textbf{R}^n)$ which is obtained by easing the restrictions on bounded
geometry in various perpendicular directions. This technical modification
allows us to inductively prove the theorem using methods of \cite{attie}.

The following definition was suggested by S.Cappell.
\begin{definition}
Let $p:X \to \textbf{R}^n$ be a control map. $p$ is said to
be $bg(r)$, if, via the decompositions $\textbf{R}^{n-r} \times \textbf{R}^r$,
$\textbf{R} \times \textbf{R}^r \times \textbf{R}^{n-r-1}$, ...,
$\textbf{R}^{n-r-1} \times \textbf{R}^r \times \textbf{R}$,
$\textbf{R}^r \times \textbf{R}^{n-r}$ the restriction
of $p$ to a regular neighborhood $D^r \times \textbf{R}^{n-r}$ of each
codimension $n-r$ hyperplane is $bg$. Thus $bg(n)$ is the same as $bg$ and
$bg(0)$ is the same as bounded control. The other notions can be considered
intermediate between $bg$ and bounded control.
\end{definition}
\begin{remark} If we consider $bg$ control to uniform control of the complexity with
complexity bound $K$, then in $\textbf{R}^2$ we can consider $bg(1)$ control to be
uniform control of complexity with complexity bound growing along lines
parallel to the $x$- or $y$-axes as either the $x$- or $y$-coordinate increases.
Note that because there must be a uniform bound on each line, the constant
is being allowed to increase along the diagonal.
In general, $bg(r)$ control is equivalent to uniform control of complexity
along hyperplanes
$\textbf{R}^{n-r}$ which is allowed to grow along an $r$-dimensional ``diagonal"
hyperplane. 
\end{remark}
\begin{definition}
$H_*^{uff}(\textbf{R}^n;Wh_*(\pi_1(M)))$ is the abelian group
$$H_0^{uff}(\textbf{R}^n;Wh(\pi_1(M))) \oplus ... \oplus 
H^{uff}_n(\textbf{R}^n;K_{1-n}(\textbf{Z}\pi_1(M)))$$
\end{definition}
\begin{theorem}
$Wh^{bg}(M \times \textbf{R}^n)=H_*^{uff}(\textbf{R}^n;Wh_*(\pi_1(M)))$
\end{theorem}
\begin{remark}
Observe that this is a ``Bass-Heller-Swan'' \cite{BHS}, \cite{fh1},
\cite{fh2},\cite{bass} formula with no Nil
terms. The reason for the absence of Nil terms in this formula is that the
``splitting'' performed in the proof below is performed with a sufficiently
large separation between hyperplanes. Almost by definition, an element of one
of the Nil terms vanishes over a large separation (i.e. larger than the
nilpotency of the element). We note further that the splitting takes place
at the level of individual elements, rather than over a set of representatives
for whole Whitehead group.
\end{remark}
Proof of 5.1. We proved this theorem for $n=1$ in \cite{attie}. We recall
the proof here. We have a map
$$Wh^{bg}(M \times \textbf{R}) \to Wh^{bdd}(M \times \textbf{R})=
\tilde{K}_0(\textbf{Z}\pi_1(M))$$
given by considering a $bg$ controlled h-cobordism as a boundedly controlled
one. We claim that if an h-cobordism is in the kernel of this map, then it
can be simultaneously split at the integer points of $\textbf{R}$. This will
be seen to follow by observing that $Wh^{bdd}(M \times \textbf{R})$ can be
identified with $\tilde{K}_0(\textbf{Z}\pi_1(M))$ which consists of splitting obstructions.

Following Pedersen, we give an explicit description of how this may be done.
Consider the integer points of $\textbf{R}$. Over each integer point is a
module $A(j)$. Suppose now we are given a controlled automorphism $\alpha$
with bound $k$. Consider the strip $l-2k \le j \le l+2k$. Define
$\phi([A,\alpha])$ by
$$\overline{A}= \bigoplus^{l+2k}_{j=l-2k}A(j)$$
$$\phi([A,\alpha])=\sum_{l \in 4k\textbf{Z}}([\overline{A},\alpha p_{-}^l
\alpha^{-1}]-[\overline{A},p_{-}^l])x_l.$$
Here $p^l_{-}$ is the projection onto the half-line below $j=l$ and $x_l$ is
the $l$-th vertex of the triangulation of $\textbf{R}$ by intervals of length
$2k$. We claim first that this defines an element of
$\tilde{K}_0(\textbf{Z}\pi_1(M))$. This follows
by inspection, since $\alpha p^l_{-}\alpha^{-1}$ differs from $p_{-}$ only
in a band around $j=l$ and these each, by inspection, can be seen to agree
for each $l$. In fact, because of the bounded geometry, one can, after uniform
stabilization, represent all of these elements by the same module. Next we
claim that this element is the splitting obstruction. We see this for the case
$l=0$, the other cases being the same. If $[A,\alpha]$ is in the kernel of the
map $\phi$ then there are modules $A^\prime$ and $A^{''}$ so that
$[\overline{A} \oplus A^{\prime} \oplus A^{''}, p_s \oplus 1 \oplus 0]$ is
isomorphic to $[\overline{A} \oplus A^{\prime} \oplus A^{''}, \alpha p_{-}
\alpha^{-1} \oplus 0]$. This implies that there is a bounded automorphism
$\beta$ so that $\beta \alpha p_{-}=p_{-}\beta \alpha$. Thus $\beta\alpha$
has the property that it preserves the module below and above $l=0$ and
hence $\beta$ preserves the half spaces below and above $l=2k+1$. This then
implies that the corresponding h-cobordism can be split by a sequence of
expansions and collapses. This follows from the fact that since one can make
the automorphism equal to the identity on each point, one can make the
h-cobordism standard over that point by a sequence of expansions and collapses.
Choosing these sufficiently far apart so that they do not interfere with
each other, we can split at each integer point.

This shows that $Wh^{bg}(M \times \textbf{R})$ can be written as a direct sum,
with one summand being represented by splitting obstructions in
$\tilde{K}_0(\textbf{Z}\pi_1(M))$ and
the kernel of the map $Wh^{bg}(M \times \textbf{R}) \to \tilde{K}_0(\textbf{Z}\pi_1(M))$ given by
h-cobordisms split over $\textbf{R}$ which is represented by elements of
$C_0^{uff}(\textbf{R};Wh(\pi_1(M)))$. We claim that two such elements are equivalent
in $Wh^{bg}$ if and only if they are homologous in $C_0^{uff}$. For if an
element is null-homologous, then there exists a bounded infinite process
trick cancelling its torsion. Conversely, if two elements are equivalent, one
can split relatively on a representative to $M \times I$ and using the new
splitting construct a homology between them. Also the map $H_0^{uff}(\textbf{R};
Wh(\pi_1(M))) \to Wh^{bg}(M \times \textbf{R})$ is given by gluing together representatives. This is
clearly injective, for if one can glue together representative h-cobordisms
so that their torsions is zero, the s-cobordism theorem shows that the
representative had to be null-homologous.

This gives a short exact sequence:
$$0 \to H_0^{uff}(\textbf{R};Wh(\pi_1(M))) \to Wh^{bg}(M \times \textbf{R}) \to
\tilde{K}_0(\textbf{Z}\pi_1(M)) \to 0$$
which is split and surjective onto the last term by the following argument.
By applying the infinite process trick of Swan, one sees that one can choose
the representatives in $\tilde{K}_0(\textbf{Z}\pi_1(M))=Wh^{bdd}(M \times \textbf{R})$ which have
bounded geometry. Given an element $(B,p) \in \tilde{K}_0(\textbf{Z}\pi_1(M))$, one constructs an
element of $Wh^{bdd}(M \times \textbf{R})$ by defining $A(j)=B$
and mapping $B$ to itself by $p$ and $1-p$ successively. This has bounded
geometry and gives rise to an automorphism whose image is $(B,p)$. Thus we
have given an element of $Wh^{bg}(M \times \textbf{R})$.

To do the general case, we first consider the case $n=2$. Observe that there 
are forgetful maps
$$Wh^{bg}(M \times \textbf{R}^2) \to Wh^{bg(1)}(M \times \textbf{R}^2)
\to Wh^{bdd}(M \times \textbf{R}^2)$$
Now note that the last group is known, due to Pedersen \cite{pedersen1},
to be $$Wh^{bdd}(M \times \textbf{R}^2) \simeq K_{-1}(\textbf{Z}\pi_1(M)).$$ 
We claim that if a geometrical representative of an
element of $Wh^{bg(1)}(M \times \textbf{R}^2)$ is the kernel of the forgetful
map $K_{-1}(M)$, then it is $bg$ split along parallel lines, and conversely.
For if an element is $bg$ split along parallel lines, then by the infinite
process trick of \cite{pedersen1}, one can move each module at each lattice
point in $Wh^{bdd}(M \times \textbf{R}^2)$ out to infinity without violating 
the loosened boundedness condition. The proof is simply that, as Pedersen shows
in \cite{pedersen1}, if a module is split along a single line in 
$Wh^{bdd}(M \times
\textbf{R}^2)$, then it is equivalent to zero, by the usual infinite process
trick, which requires no uniform bounds on the complexity of the module. 

The converse goes by an analogue of \cite{attie} with $\mathcal{C}^{bg}_1(\textbf{Z}
\pi_1(M))$ in place of $\textbf{Z}\pi_1(M)$. As in the above argument, we define the
splitting obstruction in one direction, $\phi([A,\alpha]) \in Wh^{bg(1)}(M \times
\textbf{R}^2)$ by
$$\overline{A}(j_1)=\bigoplus^{l+2k}_{j_2=l-2k} A(j_1,j_2)$$
$$\phi([A,\alpha])=\sum_{l \in 4k\textbf{Z}}([\overline{A},\alpha p^l_{-}\alpha^{-1}]-[\overline{A},p^l_{-}])x_l$$
Elements of the kernel of this map are split by the argument of \cite{pedersen1}.
Moreover the map is the forgetful map to $K_{-1}(\textbf{Z}\pi_1(M))$.
 Furthermore, applying
Swan's infinite process trick yields a surjection onto 
$K_{-1}(\textbf{Z}\pi_1(M))$ and shows
in fact that any element of $K_{-1}(\textbf{Z}\pi_1(M))$ can be represented by 
a $bg(1)$-controlled
element over $\textbf{R}^2$. Given an element 
$$(B,p) \in \tilde{K}_0^{bdd}(M \times \textbf{R}) \simeq 
K_{-1}(\textbf{Z}\pi_1(M))$$ 
one constructs an element of $Wh^{bdd}(M \times \textbf{R}^2)$ by defining
$$A(j_1,j_2)=B(j_1)$$
and mapping $B$ to itself by $p$ and $1-p$ successively, where instead of being
considered as controlled $\textbf{Z}\pi_1(M)$-modules, we consider them as  
modules over $\mathcal{C}_1(\textbf{Z}\pi_1(M))$. If we apply \cite{pedersen1}
to the category $\mathcal{C}_1(\mathcal{C}^{bg}_1(\textbf{Z}\pi_1(M)))$, we obtain for
the Whitehead group of this category $\tilde{K}_0^{bg}(M \times \textbf{R})$.
Applying the method of calculation above to this Whitehead group, we obtain
$$\tilde{K}_0^{bg}(M \times \textbf{R})=H_0^{uff}(\textbf{R};\tilde{K}_0(\textbf{Z}\pi_1(M))\oplus K_{-1}(\textbf{Z}\pi_1(M)).$$
From this we can see that $(B,p)$ can
be represented by a $bg$-controlled element over $\textbf{R}$.This has bounded
geometry in the $x$- and $y$-directions and gives rise to an automorphism 
whose image is $(B,p)$, hence it
can be thought of as an element of $Wh^{bg(1)}(M \times \textbf{R}^2)$. Thus we
have given a splitting back to  $Wh^{bg(1)}(M \times \textbf{R}^2)$.

We next analyze the map $Wh^{bg}(M \times \textbf{R}^2) \to Wh^{bg(1)}(M \times
\textbf{R}^2).$ We claim that the kernel of this map consists of elements $bg$
split simultaneously in both directions. To see this, note that if an element
is $bg(1)$ and simultaneously split, one can move the modules along the
lattice points until one reaches the diagonal and then out along the diagonals
to infinity. This is because the complexity is no longer uniformly bounded,
but is allowed to increase, remaining uniform only along parallel lines, as
one goes away from the origin. Thus the kernel of the map is represented by
elements of $C^{uff}_0(\textbf{R}^2;Wh(\pi_1(M)))$, modulo infinite process tricks
which identify different elements of $Wh^{bg}(M \times \textbf{R}^2)$, which
yields $H_0^{uff}(\textbf{R}^2;Wh(\pi_1(M)))$, as in \cite{attie}.

We will be finished with the calculation as soon as we have analyzed the kernel
of the map $Wh^{bg(1)}(M \times \textbf{R}^2) \to K_{-1}(\textbf{Z}\pi_1(M))$ and show that
both maps considered are split surjective. We claim that the obstruction to
splitting along parallel lines (separately) is given by
$$C_0^{uff}(\textbf{R};C_1^{uff}(\textbf{R};\tilde{K}_0(\textbf{Z}\pi_1(M)))\oplus C_1^{uff}(
\textbf{R};C_0^{uff}(\textbf{R};\tilde{K}_0(\textbf{Z}\pi_1(M)))$$
We clearly only have to compute this in one direction, since both
summands are isomorphic, and represent the same splitting problem. But by
the computation for the case of $\textbf{R}$, the splitting obstruction 
can be represented by an element of $C_0^{uff}(\textbf{R}; K_0(\textbf{Z}\pi_1(
M)) \oplus K_{-1}(\textbf{Z}\pi_1(M))$. The projection of this element to 
$K_{-1}(\textbf{Z} \pi_1(M))$ vanishes, by the standard infinite 
process trick, since it is split in both directions. We now observe that
since $C_1^{uff}(\textbf{R}; \textbf{Z})=\textbf{Z}$, we have the isomorphism
$$C_0^{uff}(\textbf{R};G) \simeq C_1^{uff}(\textbf{R}; C_0^{uff}(\textbf{R};G))
$$
which proves the statement.

To prove the theorem for in general for any $n$ we need to show that the kernel
of the map 
$$Wh^{bg(r)}(M \times \textbf{R}^n) \to Wh^{bg(r-1)}(M \times \textbf{R}^n)$$
is $H_{n-r}^{uff}(\textbf{R}^n;K_{1-n+r}(\textbf{Z}\pi_1(M)))$. 
We prove this inductively by splitting off an $\textbf{R}$ factor, and using
the Eilenberg-Zilber theorem. We note that the kernel of this map consists of
elements that are split in $r$ directions, since these can be taken to
infinity by an infinite process trick, along a diagonal hyperplane, 
using the loosened boundedness. The
splitting obstruction then lies in $K_0^{bg(r)}(M \times \textbf{R}^{n-1})$,
which is represented by elements of $C_0^{uff}(\textbf{R};K_0^{bg(r-1)}(
M \times \textbf{R}^{n-2}))$ by the inductive hypothesis, since the elements
of the other summands are zero, being split in more than $r$ directions. 
Using the Eilenberg-Zilber theorem completes the argument.  

\section{Proof of the Main Theorem}

We apply the splitting theorem \cite{attie}. We review the statement of this
theorem:
\begin{theorem}
Let $h:M \to N$ be a $bg$ homotopy equivalence, where $M$ is the
Cartesian product of a compact manifold with $\textbf{R}^n$ and $N$ is
$bg$ over $\textbf{R}^n$. If $X \subset N$ is a $bg$ codimension 1 submanifold
which is the transverse inverse image of $\textbf{R}^{n-1} \times \textbf{Z}$,
then $h$ can be split along $X$ if and only if an obstruction in a summand
$\tilde{K}_0^{bg}(N)$ of $Wh^{bg}(N)$ vanishes, and the components of $X$ are
sufficiently separated from each other.
\end{theorem}
\begin{definition}
We denote by $\mathcal{S}^{bg,i}_{PL}(M)$ the PL $bg$ structure set
with simple homotopy equivalences replaced by maps with torsion in
$K_i^{bg}(M)$ equal to zero. Similarly, let $\mathcal{S}^{bg,i}_{TOP}(M)$
be the topological $bg$ structure set with simple homotopy equivalences 
replaced by maps with torsion in $K_i^{bg}(M)$. 
\end{definition}
\begin{definition} Let $M$ be a compact PL manifold. Then we define the group
$H_*^{uff}(\textbf{R}^n;\mathcal{S}^{TOP}_*(M))$ to be the group
$$H_*^{uff}(\textbf{R}^n;\mathcal{S}^{TOP}_*(M))=H_0^{uff}(\textbf{R}^n;
\mathcal{S}^s_{TOP}(M \times D^n,\partial))\oplus ...$$$$...\oplus H_k^{uff}(\textbf{R}^n; \mathcal{S}^{2-k}_{TOP}(M 
\times D^{n-k}, \partial)) \oplus ... \oplus H_n^{uff}(\textbf{R}^n; 
\mathcal{S}^{2-n}_{TOP}(M))$$
Where we set the convention $\mathcal{S}^1_{TOP}(N,\partial N)=\mathcal{S}^h_{TOP}(N,\partial N)$,
for a compact PL manifold $N$ with boundary $\partial N$.
\end{definition}
\begin{definition} Let $M$ be a compact PL manifold. We define the group
$$H_*(T^n;Wh_*(G))=H_0(T^n;Wh(G)) \oplus ... \oplus H_n(T^n;K_{1-n}(\textbf{Z}G))$$
and 
$$H_*(T^n;\mathcal{S}^{TOP}_*(M))=H_0(T^n;\mathcal{S}^{s}_{TOP}(M \times D^n,\partial)) \oplus
... \oplus H_n(T^n; \mathcal{S}^{2-n}_{TOP}(M))$$
where we set the convention $\mathcal{S}^1_{TOP}(N,\partial N)=\mathcal{S}^h_{TOP}(N,\partial N)$, where $N$ is a compact PL manifold with boundary 
$\partial N$.
\end{definition}
\textit{Proof of 1.1.}We will work in the PL category and apply the PL $bg$
surgery exact sequence to get the result in the topological category.
We apply Theorem 7.1 in the following manner.

Let $N$ be simple homotopy equivalent to $M \times \textbf{R}^n$. We can
apply the splitting theorem, since the splitting obstruction vanishes, to
obtain a splitting of $N$ along parallel hyperplanes. Each split pieces
is $bg$ simple homotopy equivalent to $M \times \textbf{R}^{n-1} \times I$.
The boundary of each split piece is $bg$ simple homotopy equivalent to
$M \times \textbf{R}^{n-1}$, and each boundary piece is $bg$ PL homeomorphic
to next one. We thus have a well defined map from $\mathcal{S}^{bg,s}_{PL}(M \times
\textbf{R}^n) \to \mathcal{S}^{bg,h}_{PL}(M \times \textbf{R}^{n-1})$. A given $N$
in the kernel of this map gives rise to a chain in $C_0^{uf}(\textbf{R};
\mathcal{S}^{bg,s}_{PL}(M \times \textbf{R}^{n-1} \times I, \partial))$. By
applying the $bg$ splitting theorem, we see that two such yield $bg$ PL
homeomorphic representatives if and only if they are homologous. In addition,
one must establish the structure on the boundary of $M \times \textbf{R}^{n-1}
\times I$. 

We next perform a second splitting transverse to the first. For this we apply
the splitting theorem to each piece and then observe that the separation of
the splitting depends only on the complexity of the given homotopy equivalence,
which is uniformly bounded by hypothesis. Thus the splittings on each piece
can be aligned with each other and we obtain a transverse splitting.

We continue the analysis of the splitting as before, splitting 
$M \times \textbf{R}^{n-1} \times I$ along parallel hyperplanes, obtaining
a manifold simple homotopy equivalent to
$M \times \textbf{R}^{n-2} \times D^2$. A component of the 
 boundary of the split piece is
$bg$ homotopy equivalent to $M \times \textbf{R}^{n-2}\times I$, and the   
boundary of this manifold is $bg$ projective homotopy equivalent to $M \times
\textbf{R}^{n-2}$.
The boundary of the first splitting has been split once more, giving rise to
elements of $\mathcal{S}^{p,bg}(M \times \textbf{R}^{n-2})$ and so on.

We obtain a series of exact sequences of sets which we claim to be split:
$$0 \to C_0^{uff}(\textbf{R}; \mathcal{S}^{1-i,bg}_{PL}(M \times \textbf{R}^{n-i-1}
\times I,\partial)) \to \mathcal{S}^{1-i,bg}_{PL}(M \times \textbf{R}^{n-i})$$$$ \to
\mathcal{S}^{-i,bg}_{PL}(M \times \textbf{R}^{n-i-1})$$
the splitting being given by Cartesian product with $\textbf{R}$.
We can prove the theorem inductively by observing that 
$$\mathcal{S}^{1-i,bg}_{PL}(M \times \textbf{R}^{n-i})=C_1^{uff}(\textbf{R};
\mathcal{S}^{1-i,bg}_{PL}(M \times \textbf{R}^{n-i}))$$
and applying the PL $bg$ surgery exact sequence, the algebraic $bg$ 
surgery exact sequence, the five-lemma and the Eilenberg-Zilber theorem.

We also claim that the splitting preserves the group structures on the various
structure sets involved. This can be seen, by either using the 
``characteristic variety addition'' on the structure sets,
 or the algebraic definition of the TOP structure set. 

\textit{Proof of Theorem 1.2.} This follows from the first main theorem except for
some low dimensional difficulties, which are care of by Siebenmann periodicity, since
$$\mathcal{S}^{bg}_{TOP}(\textbf{R}^n)=H_*(\textbf{R}^n ;\mathcal{S}^{TOP}_*(pt))=0$$
where we use Siebenmann periodicity to set the convention $\mathcal{S}^{TOP}(D^i,
\partial)=0$, $i \le 4$ although the PL Poincar\'e conjecture is not known in the 
cases $i=3,4$.
(See \cite{shaneson} for a treatment of this in the compact case).
We carry this argument out following \cite{shaneson}. We then use the $bg$ PL
exact sequence and the finiteness of the group of homotopy spheres to derive
the result in the smooth category.

Let $\phi:M \to N$ be a $bg$ degree 1 normal map, $F$ a  framing of $\tau(M) \oplus 
\phi^* v$ where $\tau(M)$ is the tangent PL block bundle of $M$, $v$ the
stable normal bundle of $N$.
Let $u$ be the stable normal bundle of $CP^2$, and let $G$ be a stable framing
of $\tau(CP^2)\oplus u$. Then $F \times G$ is a framing of $$\tau(M \times
CP^2) \oplus (\phi \times id)^*(v \times u)=(\tau M \oplus \phi^*v) \times
(\tau(CP^2) \oplus u).$$ The $bg$ surgery obstruction of $M \times CP^2$, 
$\phi \times id$, $F \times G$ is then the same as that of $M$,$\phi$,$F$.
This allows us to carry out the splitting along parallel 1-, 2- and 3- dimensional
hyperplanes in $\textbf{R}^n$, by replacing $\textbf{R}^n$ with $\textbf{R}^n
\times CP^2$ and then peeling off the $CP^2$. 

\textit{Proof of Theorem 1.3.} We first observe that a map $f:M \times T^n \to N$ lifts
to a map which is boundedly homotopic to  a $bg$ simple homotopy equivalence 
on the free abelian cover if and only
if it lifts to a map which is homotopic to one which is  simple on a finite 
cover. This follows from the theorem of Bass-Heller-Swan:
$$Wh(\textbf{Z}^n \times G)=H_*(T^n;Wh_*(G))\oplus Nils$$
along with the calculation of the Whitehead group in section 5, proposition
3.4 and the observation that any given element of the Nils
vanishes on a finite cover. In fact, since $H_*^{uff}(\textbf{R}^n;G)$ is 
torsion-free for any abelian group $G$ the kernel of the map
$$H_*(T^n;Wh_*(G)) \to H_*^{uff}(\textbf{R}^n;Wh(G))$$
is torsion. This can also be seen by a transfer argument. 

We then apply the calculation of the structure set,
Theorem 1.1 above. By theorem 5.1 of \cite{shaneson}, the structure set of 
$M \times T^n$ is given by
$$\mathcal{S}^s_{TOP}(M \times T^n)=H_*(T^n;\mathcal{S}^{TOP}_*(M))$$
Applying proposition 3.4, we obtain a rational
injection  
$$H_*(T^n;\mathcal{S}_*^{TOP}(M))\otimes \textbf{Q} \to H_*^{uff}(\textbf{R}^n;\mathcal{S}_*^{TOP}(M)) \otimes \textbf{Q}.$$ 

We therefore need to check that an element of the kernel of the map, which
is a torsion element in $H_i(T^n;\mathcal{S}^{1-i}_{TOP}(M \times D^{n-i},\partial))$,
becomes trivial on a finite cover. But to see this, simply observe that taking
an $n$-fold cover of a split structure on $M \times T^n$ has the effect of
adding the structure on $(M \times D^{n-i},\partial)$ to itself $n$ times.
We can then apply the $bg$ PL surgery exact sequence.
This proves the result in the PL category. For the result in the smooth
category, observe that by \cite{hsiangshaneson, shaneson} if $N$ is PL
homeomorphic to $M \times T^n$ then a finite cover of $N$ is diffeomorphic
to $M \times T^n$.


\begin{thebibliography}{100}
\bibitem{ah}
D.Anderson and W.C.Hsiang, The functors $K_{-i}$ and pseudoisotopies of
polyhedra, Ann. of Math. 105 (1977), 201-223
\bibitem{am}
D.Anderson and H.Munkholm, Boundedly controlled topology, Lect. Notes in Math.
1323, Springer, New York (1988)
\bibitem{attie}
O.Attie, Quasi-isometry classification of some manifolds of bounded geometry,
Math. Z. 216 (1994),501-527
\bibitem{ab}
O.Attie and J.Block, Poincar\'e duality for $L^p$-cohomology, preprint.
\bibitem{abw}
O.Attie, J.Block and S.Weinberger, Characteristic classes and distortion of
diffeomorphisms, Jour. AMS 5 (1992),919-922
\bibitem{attiehurder}
O.Attie and S.Hurder, Manifolds which can not be leaves of foliations,
Topology 35 (1996) 335-353
\bibitem{bass}
H.Bass, Algebraic K-theory, W.Benjamin, New York and Amsterdam (1968)
\bibitem{BHS}
H.Bass, A.Heller and R.Swan, The Whitehead group of a polynomial extension,
Publ. Math. IHES 22 (1964), 61-79
\bibitem{baumconnes}
P.Baum and A.Connes, Geometric K-theory for Lie groups and foliations,
Preprint, IHES (1982)
\bibitem{bw}
J.Block and S.Weinberger, Aperiodic tilings, amenability and positive scalar
curvature, Jour. AMS 5 (1992),907-918
\bibitem{bw2}
J.Block and S.Weinberger, Large scale homology theories and geometry,
Geometric topology (Athens,GA,1993),522-569, AMS/IP Stud. Adv. Math. 2.1,
AMS Providence (1997)
\bibitem{browder1}
W.Browder, Torsion in H-spaces, Ann. Math. 74, 24-51 (1961)
\bibitem{browder2}
W.Browder, Homotopy type of differentiable manifolds, Proceedings Aarhus
Colloquium, 42-46 (1962)
\bibitem{cappell}
S.Cappell, A splitting theorem for manifolds and surgery groups, Invent.
Math. 33 (1976) 69-170
\bibitem{cappellw}
S.Cappell and S.Weinberger, A geometric proof of Siebenmann periodicity,
Proceedings of the 1985 Georgia Topology Conference on Geometry and
Topology, Dekker, 47-52 (1987)
\bibitem{chapman1}
T.A.Chapman, Approximation results in topological manifolds, Memoirs of the
 AMS 251
(1981)
\bibitem{chapman2}
T.A.Chapman, Controlled simple homotopy theory and applications, LNM 1009
(1983)
\bibitem{chapmanferry}
T.A.Chapman and S.C.Ferry, Approximating homotopy equivalences homeomorphisms,
Amer. J. Math. 101 (1979), 583-607
\bibitem{cheeger}
J.Cheeger, A finiteness theorem for Riemannian manifolds, Am. J. of Math.
92 (1970),62-71
\bibitem{cheegergromov}
J.Cheeger and M.Gromov, Bounds on on the von Neumann dimension of
$L^2$-cohomology and the Gauss-Bonnet theorem for open manifolds,
J. Diff Geom. 21 (1985), 1-34.
\bibitem{cms}
J.Cheeger, W. M\"uller and R.Schrader, Curvature of piecewise flat metrics,
Comm. Math. Phys. 92 (1984), 405-454
\bibitem{connellhollingsworth}
E.Connell and J.Hollingsworth, Geometric groups and Whitehead torsion,
Trans. AMS 140 (1969), 161-180
\bibitem{connesmoscovici}
A.Connes and H.Moscovici, Cyclic homology, the Novikov conjecture, and
hyperbolic groups, Topology 29 (1990), 345-388
\bibitem{dodziuk}
J.Dodziuk, de Rham-Hodge theory for $L^2$-cohomology of infinite coverings,
 Topology 16 (1977), 157-165
\bibitem{dfw}
A.N.Dranishnikov, S.Ferry, and S.Weinberger, A large Riemannian manifold
which is flexible, Ann. of Math. 157 (2003),919-938
\bibitem{fh1}
F.T.Farrell and W.C.Hsiang, A formula for $K_1 R_\alpha [T]$, Applications
of Categorical Algebra (New York, 1968), Proc. Symp. Pure Math. 17,
Amer. Math. Soc., Providence (1970), 192-218
\bibitem{fh2}
T.Farrell and W.C.Hsiang, Manifolds with $\pi_1=G \times_\alpha T.$,
Amer. J. Math. 95, 813-845 (1973)
\bibitem{ferry1}
S.Ferry, Homotoping $\epsilon$-maps to homeomorphisms, Amer. J. Math. 101
(1979), 567-582
\bibitem{fp}
S.Ferry and E.K.Pedersen, Epsilon surgery, in Novikov Conjectures, Index
Theorems and Rigidity, Oberwolfach 1993, LMS Lecture Notes Ser. 227
Ed. S.C.Ferry, A.Ranicki and J.Rosenberg (1995).
\bibitem{freedman1}
M.Freedman, The disk theorem in 4-dimensional topology, Proc. ICM, Warsaw
(1983), 647-663
\bibitem{freedman2}
M.Freedman, The topology of four-dimensional topology, J.Diff.Geometry 17
(1982), 357-453
\bibitem{fq}
M.Freedman and F.Quinn, Topology of 4-manifolds, Princeton University Press,
(1990)
\bibitem{gersten1}
S.Gersten, Bounded cohomology and group extensions, preprint (1992)
\bibitem{gersten2}
S.Gersten, Cohomological lower bounds for isoperimetric functions on groups,
Toplogy 37 (1998), 1031-1072
\bibitem{gh}
P.Griffiths and J.Harris, Principles of Algebraic Geometry, Wiley (1978)
\bibitem{gromov}
M.Gromov, K\"ahler hyperbolicity and $L^2$-Hodge theory, Jour. Diff. Geometry
33 (1991), 263-292
\bibitem{hp}
I.Hambleton and E.K.Pedersen, Bounded surgery and dihedral group actions,
Jour. of the AMS 4, 105-127 (1991)
\bibitem{hls} N.Higson, V.Lafforgue and G.Skandalis, Counterexamples to the 
Baum-Connes conjecture, GAFA 12,330-354 (2002)
\bibitem{higsonroe}
N.Higson and J.Roe, Analytic K-Homology, Oxford University Press (2000)
\bibitem{hsiangshaneson}
W.C.Hsiang and J.Shaneson, Fake tori, In: Cantrell, J. (ed.), Proceedings
of the 1969 Georgia Topology Conference, pp.18-52, Markham (1970)
\bibitem{km}
M.Kervaire and J.Milnor, Groups of homotopy spheres, Ann. of Math. 77 (1963),
504-537.
\bibitem{ks}
R.Kirby and L.Siebenmann, Essays on triangulation and smoothing of topological
manifolds, Ann. Math. Stud. (1977)
\bibitem{milnor1}
J.Milnor, On manifolds homeomorphic to the 7-sphere, Ann. of Math. 64 (1956),
399-405
\bibitem{milnor2}
J.Milnor, A procedure for killing the homotopy groups of a manifold,
Proc. Symp. Pure Math. 3 (1961), 39-55.
\bibitem{nicas}
A.Nicas, Induction theorems for groups of homotopy manifold structures,
Mem. AMS 267 (1982)
\bibitem{novikov1}
S.P.Novikov, Homotopy equivalent smooth manifolds I, Izv. Akad. Nauk SSSR
28, 365-474 (1965)
\bibitem{novikov3}
S.P.Novikov, Topological invariance of rational classes of Pontrjagin, Dokl.
Akad. Nauk SSSR 163 298-300 (1965)
\bibitem{novikov2}
S.P.Novikov, On manifolds with free abelian fundamental group applications
(Pontrjagin classes, smoothings, high-dimensional knots). Izv. Akad. Nauk
SSSR 30, 208-246 (1966)
\bibitem{pedersen1}
E.K. Pedersen, On the $K_{-i}$ functors, J.Algebra 90 (1984),461-475
\bibitem{prw}
E.K. Pedersen, J.Roe and S.Weinberger, On the homotopy invariance of the 
boundedly controlled analytic signature of a manifold over an open cone,
Novikov Conjectures, Index Theorems and Rigidity, Vol.2 (Oberwolfach 1993),
285-300 LMS Lecture Notes Ser. 227 (1995)
\bibitem{pw}
E.K. Pedersen and C.Weibel, A nonconnective delooping of algebraic K-theory,
Lect. Notes in Math. 1126 (1983)
\bibitem{quinn1}
F.Quinn, Ends of maps I, Ann. Math. 110 (1979), 275-331 
\bibitem{quinn2}
F.Quinn, Ends of maps II, Inventiones Math. 68 (1982), 353-424 
\bibitem{quinn3}
F.Quinn, Applications of topology with control, Proc. ICM (Berkeley 1986),
vol.1, Amer. Math. Soc., Providence (1987)
\bibitem{ranicki1}
A.Ranicki, Algebraic L-theory and topological manifolds, Cambridge University
Press, Cambridge (1992)
\bibitem{ranicki2}
A.Ranicki, The Hauptvermutung Book, papers by A.J.Casson, D.P.Sullivan,
M.A.Armstrong, G.E.Cooke, C.P.Rourke and A.A.Ranicki, K-Theory Book Series,
vol. 1 (1995)
\bibitem{roe1}
J.Roe, An index theorem on open manifolds I,II, Jour. Diff. Geom. 28 (1988)
87-135
\bibitem{roe2}
J.Roe, Coarse cohomology and index theory on complete Riemannian manifolds,
Mem. Amer. Math. Soc. 497 (1993)
\bibitem{rourkesanderson} 
C. Rourke and B.J. Sanderson, Block Bundles II: Transversality, Ann. of Math. 87 (1968), 256-278
\bibitem{shaneson}
J.Shaneson, Wall's surgery obstruction groups for $G \times Z$, Ann. of Math.
90 (1969)
\bibitem{siebenmann1}
L.Siebenmann, On detecting Euclidean space homotopically among topological
manifolds, Invent. Math. 6 (1968), 245-261.
\bibitem{siebenmann2}
L.Siebenmann, Topological manifolds, Proceedings 1970 Nice ICM, Gauthier-Villars,
Volume 2 (1971), 133-163 
\bibitem{sullivan}
D.Sullivan, Geometric periodicity and the invariants of manifolds, Proceedings
1970 Amsterdam Conference on Manifolds, Lecture Notes in Mathematics 197,
44-75 (1971)
\bibitem{wall}
C.T.C. Wall, Surgery on compact manifolds, Academic Press (1970)
\bibitem{weinbergeryan}
S.Weinberger and M.Yan, Equivariant periodicity for abelian group actions,
Adv. Geom. 1 (2001), 49-70 
\bibitem{whitney} 
H. Whitney, Geometric Integration Theory, Princeton University Press (1957)
\bibitem{whyte}
K.Whyte, Index theory with bounded geometry, the uniformly finite 
$\hat{A}$-class, and inifinite connected sums, Jour. Diff. Geom. 59 (2001),1-14
\end{thebibliography}
\end{document}